\documentclass[10pt]{amsart}
\usepackage{amssymb,amsmath,amsthm}
\usepackage{mathrsfs,a4wide}
\usepackage{graphicx}
\theoremstyle{plain}
\newtheorem{theorem}{Theorem}[section]
\newtheorem{lemma}[theorem]{Lemma}
\newtheorem{proposition}[theorem]{Proposition}

\newtheorem{remark}[theorem]{Remark}

\theoremstyle{definition}
\theoremstyle{remark}
\numberwithin{equation}{section}

\newcommand{\as}{{\mathcal A}}
\newcommand{\hs}{{\mathcal H}}
\newcommand{\ks}{{\mathcal K}}

\newcommand{\gs}{{\mathcal G}}
\newcommand{\fs}{{\mathcal F}}
\newcommand{\leb}{{\mathcal L}}

\newcommand{\ms}{M}

\newcommand{\eub}{{\mathcal E}}

\newcommand{\ts}{{\mathcal T}}
\newcommand{\rs}{{\mathcal R}}
\newcommand{\ws}{{\mathcal W}}

\newcommand{\tb}{{\bf T}}
\newcommand{\rb}{{\bf R}}

\newcommand{\R}{{\mathbb R}}
\newcommand{\N}{{\mathbb N}}

\newcommand{\msd}{M^{2\times2}}

\newcommand{\Om}{\Omega}
\newcommand{\Omb}{\overline{\Omega}}

\newcommand{\OmBb}{\overline{\Om}_B}

\newcommand\aplim{\mathop{\rm ap\,lim}}

\newcommand{\weak}{\rightharpoonup}

\newcommand{\eps}{\varepsilon}
\newcommand{\e}{\varepsilon}
\newcommand{\tsub}{\,\widetilde{\subset}\,}

\newcommand{\treaom}{\ts_{\varepsilon,a}(\Om)}

\newcommand{\trea}{\ts_{\varepsilon,a}(\Om')}

\newcommand{\afeaom}{\as_{\varepsilon,a}(\Om;\R^2)}

\newcommand{\afeanom}{\as^B_{\varepsilon,a_n}(\Om;\R^2)}

\newcommand{\afeaomb}{\as^B_{\varepsilon,a}(\Om;\R^2)}

\newcommand{\afenab}{\as^B_{\varepsilon_n,a}(\Om;\R^2)}
\newcommand{\afeom}{\as \fs_{\varepsilon}(\Om;\R^2)}

\newcommand{\afenabp}{\as_{\varepsilon,a}^{B}(\Om';\R^2)}

\newcommand{\Eub}{\mathcal E}
\newcommand{\Eb}{{\mathcal E}^{el}}
\newcommand{\Es}{{\mathcal E}^s}

\newcommand{\Sg}[2]{S^{#1}(#2)}

\newcommand{\radm}{\rs(\OmBb; \partial_N \Om)}
\newcommand{\gadm}{{\mathbf \Gamma}(\Om)}
\newcommand{\gadmea}{{\mathbf \Gamma}_{\eps,a}(\Om)}

\newcommand{\res}{\mathop{\hbox{\vrule height 7pt width .5pt depth 0pt
\vrule height .5pt width 6pt depth 0pt}}\nolimits}

\input psfig.sty

\title
[A discontinuous finite element approximation
]
{Discontinuous finite element approximation of \\
quasistatic crack growth in finite elasticity}
\author[A. Giacomini]
{Alessandro Giacomini}
\address[Alessandro Giacomini]{S.I.S.S.A., Via Beirut 2-4, 34014, Trieste,
Italy}
\email[A. Giacomini]{giacomin@sissa.it}
\author[M. Ponsiglione]
{Marcello Ponsiglione}
\address[Marcello Ponsiglione]{S.I.S.S.A., Via Beirut 2-4, 34014, Trieste,
Italy}
\email[M. Ponsiglione]{ponsigli@sissa.it}
\begin{document}
\vskip .2truecm
\begin{abstract}
\small{
We propose a time-space discretization of a general
notion of quasistatic growth of brittle
fractures in elastic bodies proposed in \cite{DMFT}
by G. Dal Maso, G.A. Francfort, and R. Toader,
which takes into account body forces and
surface loads. We employ adaptive triangulations and
prove convergence results for the total,
elastic and surface energies. In the case in which
the elastic energy is strictly convex,
we prove also a convergence result
for the deformations.
}
\end{abstract}
\maketitle
{\small \tableofcontents}

\section{Introduction}
\label{intr}
The aim of this paper is to provide a discontinuous finite
element approximation of a model of
quasistatic growth of
brittle fractures in finite elasticity recently
proposed in \cite{DMFT} by Dal Maso, Francfort, and Toader
in the framework of the variational theory
of crack propagation proposed by Francfort
and Marigo in \cite{FM}. This theory is inspired to Griffith's
criterion and determines the crack path through a competition
between bulk and surface energies.
\par
In the case of linearized elasticity,
a first precise mathematical formulation of the model \cite{FM}
has been given by Dal Maso and Toader \cite{DMT}:
they treat the case of
{\it anti-planar shear} in dimension two assuming that the fractures
are compact sets with a finite number of connected components.
This analysis has been extended to the case of plane elasticity
by Chambolle in \cite{Ch}. Francfort and Larsen \cite{FL},
using the framework of
$SBV$ functions (see Section \ref{notprel}),
proved the existence of a quasistatic growth of brittle
fractures in the case of {\it anti-planar shear} in any dimension
$N \ge 2$ and without
assumptions on the structure of the fractures which are dealt with
the set of jumps of the displacements.
Approximation results for the quasistatic
evolution of \cite{FL} has been given in \cite{G} and in \cite{GP}
and provide a theoretical basis to the numerical study of
the model given in \cite{BFM}.
\par
The quasistatic crack growth proposed by Dal Maso, Francfort,
and Toader in
\cite{DMFT} consider the case of finite elasticity,
and takes into account possible volume and traction forces applied
to the elastic body.
In order to describe the result of \cite{DMFT} (a complete
description is given in Section \ref{qsedmft}),
let us assume that the
elastic body has a reference configuration given by $\Om \subseteq \R^N$
open, bounded and with Lipschitz boundary. Let
$\partial_D \Om \subseteq \partial \Om$ be
open in the relative topology, and let
$\partial_N \Om := \partial \Om \setminus \partial_D \Om$.
Let $\Om_B \subseteq \Om$, and let
$\partial_S \Om \subseteq \partial_N \Om$ be such that
$\OmBb \cap \partial_S \Om=\emptyset$. $\Om_B$ is the {\it brittle part}
of $\Om$, and $\partial_S \Om$ is the part of the boundary where
traction forces are supposed to act.
A crack is given by any rectifiable set in $\OmBb$ with finite
$(N-1)$ Hausdorff measure. Given a boundary deformation $g$ on
$\partial_D \Om$ and a crack $\Gamma$,
the family of all admissible deformation of
$\Om$ is given by the set $AD(g,\Gamma)$ of all function
$u \in GSBV(\Om;\R^N)$ (see Section \ref{notprel}) such that
$S(u) \subseteq \Gamma$ and $u=g$ on $\partial_D \Om \setminus \Gamma$.
Here $S(u)$ denotes the set of jumps of $u$, and the equality $u=g$
is intended in the sense of traces. Requiring
$u=g$ only on $\partial_D \Om \setminus \Gamma$ means that
the deformation is assumed not to be transmitted through the fracture.
The bulk energy considered in \cite{DMFT} is of the form
\begin{equation*}
\int_\Om W(x,\nabla u(x))\,dx,
\end{equation*}
where $W(x,\xi)$ is quasiconvex in $\xi$, and satisfies suitable
regularity and growth assumptions (see \eqref{scisandb1}
and \eqref{scisandb2}). Moreover
the time dependent body and traction forces
are supposed to be conservative with work given by
$$
-\int_{\Om \setminus \Gamma} F(t,x,u(x))\,dx
-\int_{\partial_S \Om} G(t,x,u(x))\,d\hs^{N-1}(x),
$$
where $F$ and $G$ satisfy suitable regularity and growth
conditions (see Section \ref{qsedmft}). Finally the work
made to produce the crack $\Gamma$ is given by
$$
\Es(\Gamma):=
\int_\Gamma k(x,\nu(x))\,d\hs^{N-1}(x),
$$
where $\nu(x)$ is the normal to $\Gamma$ at $x$, and $k(x,\nu)$
satisfies
standard hypotheses which guarantee lower semicontinuity
(see Section \ref{qsedmft}). Clearly,
$W,F,G$ and $k$ depend on the material.
Let us set
$$
\Eb(t)(u):=\int_\Om W(x,\nabla u(x))\,dx
-\int_{\Om \setminus \Gamma} F(t,x,u(x))\,dx
-\int_{\partial_S \Om} G(t,x,u(x))\,d\hs^{N-1}(x),
$$
and
\begin{equation}
\label{totintr}
\Eub(t)(u,\Gamma):=\Eb(t)(u)+\Es(\Gamma).
\end{equation}
Given a boundary deformation $g(t)$ with $t \in [0,T]$ and a preexisting
crack $\Gamma_0$, a quasistatic crack growth relative to $g$ and $\Gamma_0$
is a map $\{t \to (u(t),\Gamma(t))\,:\,t \in [0,T]\}$ such that the following
conditions hold:
\begin{itemize}
\item[(1)] for all $t \in [0,T]$: $u(t) \in AD(g(t),\Gamma(t))$;
\item[]
\item[(2)] {\it irreversibility}:
$\Gamma_0 \subseteq \Gamma(s) \subseteq \Gamma(t)$
for all $0 \le s \le t \le T$;
\item[]
\item[(3)] {\it static equilibrium}: for all $t \in [0,T]$ and
for all admissible configurations $(u,\Gamma)$ with
$\Gamma(t) \subseteq \Gamma$
$$
\Eub(t)(u(t),\Gamma(t)) \le \Eub(t)(u,\Gamma);
$$
\item[]
\item[(4)] {\it nondissipativity}: the time derivative of the
total energy $\Eub(t)(u(t),\Gamma(t))$ is equal to the power of external
forces (see \eqref{nondissqse}).
\end{itemize}
\par
In this paper we discretize the model using a suitable finite
element method and prove its convergence to this notion of
quasistatic crack growth.
We restrict our analysis to a two dimensional setting considering only a
polygonal reference configuration $\Om \subseteq \R^2$.
\par
The discretization of the domain $\Om$ is carried out as in
\cite{GP} employing {\it adaptive triangulations} introduced 
by M. Negri in \cite{N} (see also \cite{N2}). Let us fix
two parameters $\eps>0$ and $a \in ]0,\frac{1}{2}[$.
We consider a regular triangulation $\rb_\eps$ of size $\eps$
of $\Om$, i.e. we assume that there exist two constants
$c_1$ and $c_2$ so that
every triangle $T \in \rb_\eps$ contains a ball of diameter
$c_1 \eps$ and is contained in a ball of diameter $c_2 \eps$.
In order to treat the boundary data,
we assume also that $\partial_D \Om$ is composed of edges of $\rb_\eps$.
On each edge $[x,y]$ of $\rb_\eps$ we consider a point $z$
such that $z=tx+(1-t)y$ with
$t \in [a, 1-a]$. These points are called {\it adaptive vertices}.
Connecting
together the adaptive vertices, we divide every $T \in \rb_\eps$
into four triangles.
We take the new triangulation $\tb$ obtained after this
division as the discretization
of $\Om$. The family of all such triangulations will be denoted by
$\ts_{\eps,a}(\Om)$.
\par
The discretization of the energy functional is obtained restricting
the total energy \eqref{totintr} to the family of functions $u$
which are affine on the triangles of some triangulation
$\tb(u) \in \ts_{\eps,a}(\Om)$
and are allowed to jump across the edges of $\tb(u)$ contained
in $\OmBb$.
We indicate this space by $\afeaomb$.
The boundary data is assumed to belong
to the space $\afeom$ of continuous functions which are affine
on every triangle $T \in
\rb_\eps$.
\par
Given the boundary data
$g_\eps \in W^{1,1}([0,T], W^{1,p}(\Om;\R^2) \cap L^q(\Om;\R^2))$
with $g_\eps(t) \in \afeom$ for all $t \in [0,T]$ ($p,q$ are related
to the growth assumptions on $W,F,G$) and an initial crack
$\Gamma^0_{\eps,a}$ (see Section \ref{devol}),
we divide $[0,1]$ into subintervals $[t^\delta_i,t^\delta_{i+1}]$
of size $\delta>0$ for $i=0, \ldots, N_\delta$, and
for all $u \in \afeaomb$ we indicate by $S_D^{g_\eps(t)}(u)$ the edges
of the triangulation $\tb(u)$ contained in
$\partial_D \Om$ on which $u \not= g_\eps(t)$.
Using a variational argument (Proposition \ref{discrevol}),
we construct a {\it discrete evolution}
$\{(u^{\delta,i}_{\eps,a},\Gamma^{\delta,i}_{\eps,a})\,:
\,i=0, \ldots,N_\delta\}$ such that for all $i=0, \ldots, N_\delta$
we have
$u^{\delta,i}_{\eps,a} \in \afeaomb$,
\begin{equation*}
\Gamma^{\delta,i}_{\eps,a}:=
\bigcup_{r=0}^i \big[ S(u^{\delta,r}_{\eps,a}) \cup
S_D^{g_\eps(t^\delta_r)}(u^{\delta,r}_{\eps,a}) \big],
\end{equation*}
and the following {\it unilateral minimality property} holds:
for all $v \in \afeaomb$
\begin{equation}
\label{pieceminintr}
\Eb(t^\delta_i)(u^{\delta,i}_{\eps,a}) \le
\Eb(t^\delta_i)(v)+
\Es\left( \big( S(v) \cup S_D^{g_\eps(t^\delta_i)}(v) \big) \setminus
\Gamma^{\delta,i-1}_{\eps,a}\right).
\end{equation}
Notice that by construction $u^{\delta,i}_{\eps,a} \in
AD(g_\eps(t^\delta_i),\Gamma^{\delta,i}_{\eps,a})$.
Moreover the definition
of the discrete fracture ensures that
$\Gamma^{\delta,i}_{\eps,a} \subseteq
\Gamma^{\delta,j}_{\eps,a}$
for all $i \le j$, recovering in this discrete
setting the irreversibilty
of the crack growth given in $(2)$.
The minimality property \eqref{pieceminintr}
is the reformulation in the finite element space of
the equilibrium condition $(3)$.
Finally we obtain an estimate from above for
$\Eub(t^\delta_i)(u^{\delta,i}_{\eps,a},\Gamma^{\delta,i}_{\eps,a})$
(see Proposition \ref{dener}) which is a discrete version of $(4)$.
\par
In order to perform the asymptotic analysis of the
{\it discrete evolution}
$\{(u^{\delta,i}_{\eps,a},\Gamma^{\delta,i}_{\eps,a})
\,:\,i=0, \ldots,N_\delta\}$ we make
the piecewise constant interpolation in time
$u^{\delta}_{\eps,a}(t)=u^{\delta,i}_{\eps,a}$ and
$\Gamma^{\delta}_{\eps,a}(t)=\Gamma^{\delta,i}_{\eps,a}$ for all
$t^\delta_i \le t <t^\delta_{i+1}$.
Let us suppose that
$$
g_\eps \to g
\quad
\quad
\text{strongly in }W^{1,1}([0,T],W^{1,p}(\Om;\R^2) \cap L^q(\Om;\R^2))
$$
(where on $W^{1,p}(\Om;\R^2) \cap L^q(\Om;\R^2)$ we take the norm
$\|u\|:=\|u\|_{W^{1,p}(\Om;\R^2)}+\|u\|_{L^q(\Om;\R^2)}$),
and that $\Gamma^0_{\eps,a}$ approximate an initial crack
$\Gamma^0$ in the sense of Proposition \ref{gammazero}.
\par
The main result of the paper (Theorem \ref{mainthm}) states that
there exist a quasistatic evolution $\{t \to (u(t),\Gamma(t))\,:\,
t \in [0,T]\}$
in the sense of \cite{DMFT} relative to the boundary deformation $g$
and the preexisting crack $\Gamma^0$
and sequences $\delta_n \to 0$, $\eps_n \to 0$,
$a_n \to 0$,
such that setting
$$
u_n(t):=u^{\delta_n}_{\eps_n,a_n}(t),
\quad\quad
\Gamma_n(t):=\Gamma^{\delta_n}_{\eps_n,a_n}(t),
$$
for all $t \in [0,T]$ the following facts hold:
\begin{itemize}
\item[]
\item[(a)] $(u_n(t))_{n \in \N}$
is weakly precompact in
$GSBV^p_q(\Om;\R^2)$, and every
accumulation point $\tilde{u}(t)$ is such that
$\tilde{u}(t) \in AD(g(t),\Gamma(t))$, and $(\tilde{u}(t),\Gamma(t))$
satisfy the static equilibrium $(2)$; moreover there exists a 
subsequence $(\delta_{n_k}, \eps_{n_k}, a_{n_k})_{k \in \N}$
of $(\delta_n, \eps_n, a_n)_{n \in \N}$ (depending on $t$)
such that
$$
u_{n_k}(t) \weak u(t)
\quad
\text{ weakly in } GSBV^p_q(\Om; \R^2)
$$
(see Section \ref{notprel} for a precise definition of
$GSBV^p_q(\Om;\R^2)$, and of weak convergence in this space);
\item[]
\item[(b)]
convergence of the total energy holds, and
more precisely elastic and surface energies
converge separately, that is
\begin{equation*}
\Eb(t)(u_n(t)) \to \Eb(t)(u(t)),
\quad\quad
\Es(\Gamma_n(t)) \to \Es(\Gamma(t)).
\end{equation*}
\end{itemize}
By point $(a)$, the approximation of the deformation $u(t)$ is available only
up to a subsequence depending on $t$: this is due to the possible non uniqueness
of the minimum energy deformation associated to $\Gamma(t)$. In the case
$\Eb(t)(u)$ is strictly convex, it turns out that
the deformation $u(t)$ is uniquely determined, and we prove that
(Theorem \ref{main2})
$$
\nabla u_n(t) \to \nabla u(t)
\quad\quad
\text{strongly in }L^p(\Om;\msd),
$$
and
$$
u_n(t) \to u(t)
\quad\quad
\text{strongly in }L^q(\Om;\R^2).
$$
The main difficulty to prove Theorem \ref{mainthm} consists
in passing to the limit in the static equilibrium
\eqref{pieceminintr}. In order to find the fracture $\Gamma(t)$
in the limit, in Lemma \ref{agamma} and Lemma \ref{ato0}
we adapt to the context of finite elements
the notion of $\sigma^p${-} convergence of sets formulated in \cite{DMFT}.
This is the key tool to obtain the convergence of elastic and
surface energies at all times $t \in [0,T]$ (while
in \cite{GP} this was available only at the continuity points
of $\hs^1(\Gamma(t))$).
In order to infer the static equilibrium of $\Gamma(t)$
from that of $\Gamma_n(t)$, we employ a generalization of the
piecewise affine transfer of jumps \cite[Proposition 5.1]{GP}
(see Proposition \ref{piecetransf2}).
\par
The paper is organized as follows. In Section \ref{notprel}
we introduce the basic notation, and some tools employed
throughout the paper. In Section \ref{qsedmft}
we describe the quasistatic crack growth of \cite{DMFT}
precising the functional setting and the hypotheses on
the elastic and surface energies involved. In Section
\ref{femspace} we introduce the finite element space, and
in Section \ref{secincrack} we prove an approximation result for
a preexisting crack configuration.
In Section \ref{devol} we prove the existence of a discrete
evolution, and in Section \ref{approx} we prove the main
approximation result (Theorem \ref{mainthm}). In Section
\ref{convexcase} we treat the case of strictly convex total energy.

\section{Notations and Preliminaries}
\label{notprel}

In this section we introduce the main notations and the
preliminary results employed in
the rest of the paper.

\vskip10pt\noindent
{\bf Basic notation.}
We will employ the following basic notation:
\begin{itemize}
\item[-] $\ms^{n\times m}$ is the space of
${n\times m}$ matrices;
\item[-] $\hs^{1}$ is the one-dimensional Hausdorff measure;
\item[-] for $p \in [1,+\infty]$, 
$\|\cdot\|_p$ denotes the usual $L^p$ norm;
\item[-] if $\mu$ is a measure on $\R^2$ and $A$ is a
Borel subset of $\R^2$,
$\mu \res A$ denotes the restriction of $\mu$ to $A$, i.e.
$(\mu \res A)(B):=\mu(B \cap A)$ for all Borel sets
$B \subseteq \R^2$;
\item[-] if $A,B \subseteq \R^2$, $A \tsub B$ means that
$A \subseteq B$ up to a set of $\hs^1${-}measure zero.
\end{itemize}

\vskip20pt\noindent
{\bf SBV and GSBV spaces.}
Let $A$ be an open subset of $\R^n$, and let
$u: A \to \R^m$ be a measurable function.
Given $x \in A$, we say that $\tilde{u}(x)$ is the
{\it approximate limit}
of $u$ at $x$, and we write 
$\tilde{u}(x)=\aplim\limits_{y \to x}u(y)$, 
if for every $\eps>0$
$$
\lim_{r \to 0} r^{-n} \leb^n 
\left( \{y \in B_r(x)\,:\,
|u(y)-\tilde{u}(x)|>\eps\} \right)=0.
$$
Here $B_r(x)$ denotes the ball of center $x$ and radius $r$.
We indicate by $S(u)$ the set of points where the approximate
limit of $u$ does not exist.
We say that the matrix $m \times n$ $\nabla u(x)$ is the
approximate gradient of $u$ at $x$ if
$$
\aplim\limits_{y \to x}
\frac{u(y)-u(x)-\nabla u(x)(y-x)}{|y-x|}=0.
$$
\par
We say that $u \in BV(A;\R^m)$ if $u \in L^1(A;\R^m)$,
and its distributional
derivative $Du$ is a vector-valued Radon measure on $A$.
In this case, it turns out
that $S(u)$ is rectifiable, that is there exists a sequence
$(M_i)_{i \in \N}$ of $C^1${-}manifolds such that
$S(u) \subseteq \bigcup_i M_i$ up to a set of
$\hs^{n-1}${-}measure zero;
as a consequence $S(u)$ admits a normal $\nu_x$ for
$\hs^{n-1}${-}almost every
$x \in S(u)$. Moreover the approximate gradient
$\nabla u(x)$ exists for a.e. $x \in A$, and
$\nabla u$ is the density of the absolutely
continuous part of $Du$.
\par
We say that $u \in SBV(A;\R^m)$ if $u \in BV(A;\R^m)$
and the singular
part $D^su$ of its distributional derivative $Du$ is
concentrated on $S(u)$.
The space $SBV(A;\R^m)$ is called the space of $\R^m$-valued
{\it special functions of bounded variation}.
For more details, the reader is referred to \cite{AFP}.
We indicate with $SBV_{loc}(A,\R^m)$ the space of functions
which belong to $SBV(A',\R^m)$
for every open set $A'$ with compact closure in $A$.
\par
The set $GSBV(A,\R^m)$ is defined as the set of functions
$u:A \to \R^m$ such
that $\varphi(u) \in SBV_{loc}(A)$ for every
$\varphi \in C^1(\R^m)$
such that the support of $\nabla \varphi$ has compact
closure in $\R^m$.
If $p \in ]1,+\infty[$, we set
\begin{equation*}
GSBV^p(A,\R^m):=
\{u\in GSBV(A,\R^m)\,:\, \nabla u \in L^p(A,\ms^{m\times n}),\,
\hs^{n-1}(S(u)) <+\infty\}.
\end{equation*}
By \cite[Proposition 2.2]{DMFT} the space $GSBV^p(A,\R^m)$
coincide with
$(GSBV^p(A,\R))^m$, that is $u:=(u_1, \dots,u_m) \in
GSBV^p(A,\R^m)$
if and only if $u_i \in GSBV^p(A,\R)$ for every
$i=1,\dots,m$.
\par
The following compactness and lower semicontinuity result
will be used in the following sections.
For a proof, we refer to \cite{A2}.

\begin{theorem}
\label{GSBVcompact}
Let $A$ be an open and bounded subset of $\R^n$. Let
$g(x,u): \, A\times \R^m \to [0,\infty]$ be a Borel function,
lower semicontinuous in $u$ and satisfying the condition
$$
\lim_{|u|\to \infty} g(x,u)=+\infty \text{ for a.e. } x \in A.
$$
Let $(u_k)_{k \in \N}$ be a sequence in $GSBV^p(A;\R^m)$ such that
$$
\limsup_k \int_A |\nabla u_k(x)|^p \,dx+
\hs^{n-1} \left( S(u_k) \right)+
\int_A g(x, u_k(x)) \, dx <+\infty.
$$
Then there exists a subsequence $(u_{k_h})_{h \in \N}$ and
a function $u \in GSBV^p (A;\R^m)$ such that
\begin{align}
\label{gsbvconv}
u_{k_h} \to u &\quad\quad\text{in measure}, \\ \nonumber
\nabla u_{k_h} \weak \nabla u &\quad\quad\text{weakly in }
L^p(A;M^{m \times n}).
\end{align}
Moreover we have that
\begin{equation*}
\hs^{n-1} \left( S(u) \right)
\le \liminf_h
\hs^{n-1} \left( S(u_{k_h}) \right).
\end{equation*}
\end{theorem}

Let $q \in ]1,+\infty[$ and let us set
\begin{equation}
\label{defgsbvpq}
GSBV^p_q(A;\R^m):=GSBV^p(A;\R^m) \cap L^q(A;\R^m).
\end{equation}
We say that $u_k \weak u$ weakly in $GSBV^p_q(A;\R^m)$ if
\begin{align}
\nonumber
u_k \to u &\quad\quad\text{in measure} \\
\nabla u_k \weak \nabla u &\quad\quad\text{weakly in }
L^p(A; M^{m \times n}) \\
\nonumber
u_k \weak u &\quad\quad\text{weakly in }L^q(A;\R^m).
\end{align}
We will often use the following fact: if $u_k \weak u$
weakly in $GSBV^p_q(A;\R^m)$ and $\Gamma \subseteq A$
is such that $\hs^{N-1}(\Gamma)<+\infty$ and $S(u_k) \subseteq \Gamma$
up to a set of $\hs^{N-1}${-}measure zero for all $k$,
then $S(u) \subseteq \Gamma$ up to a set of $\hs^{N-1}${-}measure zero.

\vskip20pt\noindent
{\bf $\Gamma$-convergence.}
Let us recall the definition of De Giorgi's
{\it $\Gamma$-convergence} in metric spaces:
we refer the reader to \cite{dm} for an exhaustive
treatment of this subject.
Let $(X,d)$ be a metric space. We say that a sequence
$F_h:X\to [-\infty ,+\infty
]$ $\Gamma $-converges to $F:X\to [-\infty ,+\infty ]$
(as $h\to +\infty$) if for all $u \in X$ we have
\begin{itemize}
\item[{\rm (i)}] ({\it $\Gamma${-}liminf inequality})
for every sequence $(u_h)_{h \in \N}$ converging to
$u$ in $X$,
$$
\liminf\limits _{h\to+\infty }F_h(u_h)\geq F(u);
$$
\item[{\rm (ii)}] ({\it $\Gamma${-}limsup inequality})
there exists a sequence
$(u_h)_{h \in \N}$ converging to $u$ in $X$, such that
$$
\limsup\limits _{h\to +\infty }F_h(u_h)\leq F(u).
$$
\end{itemize}
The function $F$ is called the $\Gamma${-}limit of $(F_h)$
(with respect to $d$),
and we write $F\,=\,\Gamma{-}\lim_{h}F_h$.
\par
We say that a family of functionals $\{F_\eps\}$
$\Gamma${-}converges to $F$ as
$\eps \to 0$ if for every sequence $\eps_h \to 0$ as
$h \to +\infty$ we have
$\Gamma{-}\lim_h F_{\eps_h}=F$.
\par
The peculiarity of this type of convergence
is its variational
character explained in the following proposition.

\begin{proposition}
Assume that the sequence $(F_h)_{h \in \N}$
$\Gamma${-}converges to
$F$ and that there
exists a compact set $K\subseteq X$ such that for all
$h \in \N$
$$
\inf\limits _{u\in K}F_h(u)=\inf\limits _{u\in X}
F_h(u).
$$
Then $F$ admits a minimum on $X$,
$\inf_{X}F_h \to \min_X F$, and
any limit point of any sequence $(u_h)_{h \in \N}$ such that
$$
\lim\limits _{h\to +\infty }\Bigl( F_h(u_h)-
\inf\limits _{u\in X}F_h(u)\Bigr) =0
$$
is a minimizer of $F$.
\end{proposition}

\vskip20pt\noindent
{\bf Hausdorff metric on compact sets.}
Let $A \subseteq \R^2$ be open and bounded, and
let $\ks(\overline{A})$
be the set of all compact subsets of $\overline{A}$.
$\ks(\overline{A})$ can be endowed by the
Hausdorff metric $d_H$ defined by
$$
d_H(K_1,K_2) :=
\max \left\{ \sup_{x \in K_1} {\rm dist}(x,K_2),
\sup_{y \in K_2} {\rm dist}(y,K_1)\right\},
$$
with the conventions
${\rm dist}(x, \emptyset)= {\rm diam}(A)$ and
$\sup\emptyset=0$, so that $d_H(\emptyset, K)=0$ if
$K=\emptyset$ and
$d_H(\emptyset,K)={\rm diam}(A)$ if $K \not=\emptyset$.
It turns out that
$\ks(\overline{A})$ endowed with the Hausdorff metric
is a compact space
(see e.g. \cite{Ro}).

\section{The quasistatic crack growth of
Dal Maso-Francfort-Toader}
\label{qsedmft}
In this section we describe the quasistatic evolution of brittle
fractures proposed in
\cite{DMFT}. They consider the case of  $n$-dimensional
finite elasticity, for an
arbitrary $n\ge1$, with a quasiconvex bulk energy and with
prescribed boundary deformations and applied loads, depending on
time. Since we are going to approximate the case $n=2$,
we prefer to introduce the model
in this particular case. For more details, we refer the
reader to \cite{DMFT}.
\par
Let $\Om$ be a bounded open set of $\R^2$ with Lipschitz boundary
and let $\Om_B$
be an open subset of $\Om$. Let $\partial_N \Om \subseteq \partial \Om$
be closed in the relative topology,
and let $\partial_D\Om:= \partial\Om\setminus \partial_N\Om$.
Let $\partial_S\Om \subseteq \partial_N\Om$ be closed in the
relative topology and such that
$\OmBb\cap \partial_S\Om = \emptyset$.
In the model proposed in \cite{DMFT},
$\Om_B$ represents the brittle part of $\Om$, and $\partial_D \Om$
the part of the boundary on which
the deformation is prescribed. Moreover the elastic body $\Om$ is
supposed to be subject to surface forces
acting on $\partial_S \Om$.
\par\medskip\noindent
{\bf Admissible cracks and deformations.}
The set of admissible cracks is given by
\begin{equation*}
\radm:= \{\Gamma:\Gamma\text{ is rectifiable },
\Gamma \tsub (\OmBb \setminus \partial_N\Om),
\, \hs^1(\Gamma)<+\infty\}.
\end{equation*}
Here $A \tsub B$ means that $A \subseteq B$ up to a set of
$\hs^1$-measure zero, and $\Gamma$ rectifiable
means that there exists a sequence $(M_i)$ of $C^1${-}manifolds 
such that $\Gamma \tsub \bigcup_i M_i$.
If $\Gamma$ is rectifiable, we can define normal vector fields 
$\nu$ to $\Gamma$ in the following way:
if $\Gamma=\bigcup_i \Gamma_i$ with $\Gamma_i \tsub M_i$ and 
$\Gamma_i \cap \Gamma_j=\emptyset$ for $i \not= j$,
given $x \in \Gamma_i$, we take $\nu(x)=\nu_{M_i}(x)$, where
$\nu_{M_i}(x)$ is a normal vector to the
$C^1${-}manifold $M_i$ at $x$. It turns out that 
two normal vector fields
associated to different decompositions 
$\bigcup_i \Gamma_i$ of $\Gamma$
coincide up to the sign $\hs^{1}$ almost everywhere.
\par
Given a crack $\Gamma$, an admissible deformation is given by 
any function $u \in
GSBV(\Om;\R^2)$ such that $S(u) \tsub \Gamma$.
\par\medskip\noindent
{\bf The surface energy.}
The surface energy of a crack $\Gamma$ is given by
\begin{equation}
\label{crackener}
\Es(\Gamma):=\int_\Gamma k(x,\nu(x)) \,d\hs^1(x),
\end{equation}
where $\nu$ is a unit normal vector field on $\Gamma$.
Here $k:\OmBb\times\R^2 \to \R$ is continuous, 
$k(x,\cdot)$ is a norm in $\R^2$
for all $x \in \OmBb$ and for all
$x\in  \OmBb$
and $\nu\in \R^2$
\begin{equation}
\label{ik3}
K_1|\nu| \le k(x, \nu) \le K_2|\nu|,
\end{equation}
where $K_1,K_2>0$. Notice that since $k$ is even in the
second variable,
we have that the integral \eqref{crackener} is independent 
of the
orientation given to $\Gamma$, that is independent of 
the particular choice of the unit
normal vector field $\nu$.
\par\medskip\noindent
{\bf The bulk energy.}
Let $p>1$ be fixed.
Given a deformation $u \in GSBV^p(\Om;\R^2)$ the
associated {\it bulk energy}
is given by
\begin{equation}
\label{bws}
\ws(\nabla u):= \int_\Om W(x,\nabla u(x)) \, dx,
\end{equation}
where $W: \Om \times \msd \to [0,+\infty)$ is a
Carath\'eodory function satisfying
\begin{align}
\label{scisandb1}
& \text{ for every } x \in \Om:\,
W(x,\cdot) \text{ is quasiconvex and }
C^1 \text{ on }\msd,\\
\label{scisandb2}
& \text{ for every } (x,\xi) \in \Om\times \msd:
a_0^W|\xi|^p - b_0^W(x)\le W(x,\xi)\le
a_1^W |\xi|^p + b_1^W(x).
\end{align}
Here $a_0^W, a_1^W>0$, and
$b_0^W, b_1^W \in L^1(\Om)$
are nonnegative functions. Quasiconvexity
of $W$ means that for all $\xi \in \msd$ and
for all $\varphi \in C^\infty_c(\Om;\R^2)$
$$
W(\xi) \le \int_\Om W(\xi +\nabla \varphi)\,dx.
$$
If we denote by
$\partial_\xi W: \Om \times \msd \to \msd$
the partial derivative of $W$ with respect to
$\xi$, it turns out that
there exists a positive constant $a_2^W>0$ and a
nonnegative function $b_2^W\in L^{p'}(\Om)$, with
$p':=p/(p-1)$, such that for all
 $(x,\xi)\in \Om\times  \msd$
\begin{equation}
\label{estgradxi}
|\partial_\xi W(x,\xi)|\le a_2^W |\xi|^{p-1} + b_2^W(x).
\end{equation}
By \eqref{scisandb2} and  \eqref{estgradxi} the functional
$\ws$, defined for all $\Phi \in L^p(\Om; \msd)$ by
$$
\ws(\Phi):=\int_\Om W(x,\Phi(x)) \, dx,
$$
is of class $C^1$ on $L^p(\Om;\msd)$, and its differential
$\partial \ws:L^p(\Om;\msd) \to L^{p'}(\Om;\msd)$ is given by
\begin{equation*}
\langle \partial \ws(\Phi),\Psi \rangle =
\int_\Om \partial_\xi W(x,\Phi(x)) \Psi(x)\, dx,
\quad\quad
\Phi,\Psi \in L^p(\Om;\msd),
\end{equation*}
where $\langle \cdot,\cdot \rangle$ denotes the duality
pairing between the spaces
$L^{p'}(\Om;\msd)$ and $L^p(\Om;\msd)$.
By \eqref{scisandb2} and  \eqref{estgradxi},
there exist six positive constants
$\alpha_0^{\ws}>0$, $\alpha_1^{\ws}>0$,
$\alpha_2^{\ws}>0$, $\beta_0^{\ws} \ge 0$,
$\beta_1^{\ws} \ge 0$,
$\beta_2^{\ws} \ge 0$ such that for every
$\Phi,\,\Psi\in L^p(\Om;\msd)$
\begin{equation*}
\alpha_0^{\ws}\|\Phi\|_p^p - \beta_0^{\ws} \le
\ws(\Phi)\le \alpha_1^{\ws}\|\Phi\|_p^p
+\beta_1^{\ws},
\end{equation*}
\begin{equation}
\label{wsb}
|\langle \partial \ws(\Phi),\Psi \rangle|\le
(\alpha_2^{\ws}\|\Phi\|_p^{p-1} + \beta_2^{\ws})
\|\Psi\|_p.
\end{equation}
\par\medskip\noindent
{\bf The body forces.} Let $q>1$ be fixed.
The density of applied
body forces per unit volume in the
reference configuration relative to the deformation
$u$ at time $t \in [0,T]$
is given by $\partial_z F(t,x,u(x))$.
Here $F:[0,T]\times \Om \times \R^2 \to \R$ is such that:
\begin{align*}
&\text{ for every } z\in \R^2:\, (t,x)\to F(t,x,z)
\text{ is }\leb^1 \times \leb^2 \text{ measurable on }
[0,T]\times \Om, \\
&\text{ for every } (t,x)\in [0,T] \times \Om:\,
z \to F(t,x,z) \text{ belongs to } C^1(\R^2),
\end{align*}
and satisfies the following growth conditions
\begin{align}
\label{hpf4}
&a_0^F|z|^q - b_0^F(t,x) \le -F(t,x,z)\le a_1^F|z|^q + b_1^F(t,x), \\
\nonumber
&|\partial_z F(t,x,z)| \le a_2^F |z|^{q-1} + b_2^F(t,x)
\end{align}
for every $(t,x,z) \in [0,T]\times \Om \times \R^2$, with
$a_0^F>0$, $a_1^F>0$ and $a_2^F>0$, and
where $b_0^F, b_1^F \in C^0([0,T];L^1(\Om))$,
$b_2^F \in C^0([0,T];L^{q'}(\Om))$
are nonnegative functions, with $q':=q/(q-1)$.
\par
In order to deal with time variations, we assume also that
for every $(t,z) \in [0,T] \times \R^2$
\begin{align*}
&F(t,x,z) = F(0,x,z)+ \int_0^t \dot{F}(s,x,z) \,ds
\quad\quad \text{ for a.e. }
x \in \Om, \\
\nonumber
&\partial_z F(t,x,z)= \partial_z F(0,x,z)+
\int_0^t \partial_z \dot{F}(s,x,z) \,ds \quad\quad
\text{ for a.e. }x \in \Om,
\end{align*}
where $\dot{F}:[0,T] \times \Om \times \R^2 \to \R$ is
such that
\begin{align*}
& \text{ for all } z \in \R^2:\, (t,x) \to \dot{F}(t,x,z)
\text{ is }\leb^1 \times \leb^2
\text{ measurable on } [0,T]\times \Om,  \\
& \text{ for all }(t,x) \in [0,T] \times \Om:\,
z \to \dot{F}(t,x,z)
\text{ is of class $C^1$ on }\R^2,
\end{align*}
and satisfies the growth conditions
\begin{align*}
& |\dot{F}(t,x,z)|\le a_3^F(t)|z|^{\dot{q}} + b_3^F(t,x), \\
& |\partial_z \dot{F}(t,x,z)|\le a_4^F(t)|z|^{\dot{q}-1}
+ b_4^F(t,x)
\end{align*}
for all $(t,x,z) \in [0,T] \times \Om \times \R^2$. Here
$1 \le \dot{q} < q$, and $a_3^F, a_4^F \in L^1([0,T])$,
$b_3^F \in L^1([0,T];L^1(\Om))$,
$b_4^F \in L^1([0,T]; L^{\dot{q}'}(\Om))$ are nonnegative
functions with
$\dot{q}':=\frac{\dot{q}}{\dot{q}-1}$.
\par
Under the previous assumptions, for every $t\in [0,T]$
the functionals
\begin{equation}
\label{bodyenergy}
\fs(t)(u):= \int_\Om F(t,x,u(x))\,dx, \quad\quad
\dot{\fs}(t)(u):= \int_\Om \dot{F}(t,x,u(x))\,dx
\end{equation}
are well defined on $L^q(\Om;\R^2)$ and $L^{\dot{q}}(\Om;\R^2)$
respectively.
Moreover we have that $\fs(t)$ is of class $C^1$ on
$L^q(\Om;\R^2)$, with differential $\partial \fs(t):
L^q(\Om;\R^2)\to L^{q'}(\Om;\R^2)$
defined by
$$
\langle \partial \fs(t)(u),v \rangle =
\int_\Om \partial_z F(t,x,u(x)) v(x)\, dx,
\quad\quad
u,v \in L^q(\Om;\R^2),
$$
where $\langle \cdot,\cdot \rangle$ denotes now the duality
pairing between
$L^{q'}(\Om;\R^2)$ and $L^q(\Om;\R^2)$.
$\dot{\fs(t)}$ is $C^1$ on $L^{\dot{q}}(\Om;\R^2)$ with
differential defined by
$$
\langle \partial \dot{\fs}(t)(u),v \rangle =
\int_\Om \partial_z \dot{F}(t,x,u(x)) v(x)\, dx,
\quad\quad
u,v \in L^{\dot{q}}(\Om;\R^2),
$$
where $\langle \cdot,\cdot \rangle$ denotes the duality
pairing between
$L^{\dot{q}'}(\Om;\R^2)$ and $L^{\dot{q}}(\Om;\R^2)$.
For all $u,v \in L^q(\Om;\R^2)$ and for all $t \in [0,T]$
we have
\begin{equation*}
\fs(t)(u) = \fs(0)(u)+ \int_0^t \dot{\fs}(s)(u) \,ds,
\end{equation*}
\begin{equation}
\label{ac2bis}
\langle \partial \fs(t)(u),v\rangle=
\langle \partial \fs(0)(u),v\rangle+
\int_0^t \langle\partial \dot{\fs}(s)(u),v\rangle \,ds.
\end{equation}
Moreover we have that for every $t\in [0,T]$ and for every $u$,
$v\in L^q(\Om;\R^n)$
\begin{align}
\nonumber
& \alpha_0^{\fs} \|u\|_q^q - \beta_0^{\fs} \le -\fs(t)(u)
\le \alpha_1^{\fs}\|u\|_q^q + \beta_1^{\fs},\\
\label{fs2}
& |\langle\partial \fs(t)(u),v\rangle|\le
(\alpha_2^{\fs}\|u\|_q^{q-1}
+\beta_2^{\fs})\|v\|_q,\\
\label{fs3}
& |\dot{\fs}(t)(u)|\le \alpha_3^{\fs}(t)
\|u\|_{\dot{q}}^{\dot{q}} + \beta_3^{\fs}(t), \\
\label{fs4}
& |\langle\partial \dot{\fs}(t)(u),v\rangle|\le
(\alpha_4^{\fs}(t) \|u\|_{\dot{q}}^{\dot{q}-1}+\beta_4^{\fs}(t))
\|v\|_{\dot{q}} ,
\end{align}
where $\alpha_0^{\fs}>0$, $\alpha_1^{\fs}>0$, $\alpha_2^{\fs}>0$,
$\beta_0^{\fs} \ge 0$, $\beta_1^{\fs} \ge 0$,
$\beta_2^{\fs} \ge 0$ are positive constants,
and $\alpha_3^{\fs}, \alpha_4^{\fs}, \beta_3^{\fs},
\beta_4^{\fs} \in L^1([0,T])$
are nonnegative functions.
\par\medskip\noindent
{\bf The surface forces.}
The density of the surface forces on $\partial_S \Om$ at
time $t$ under the
deformation $u$ is given by $\partial_z G(t,x,u(x))$,
where $G:[0,T] \times \partial_S \Om \times \R^2 \to \R$
is such that
\begin{align*}
& \text{ for every } z \in \R^2:\,(t,x) \to G(t,x,z)
\text{ is } \leb^1 \times \hs^1\text{{-}measurable},\\
& \text{ for every } (t,x) \in [0,T] \times \partial_S\Om:\,
z \to G(t,x,z) \text{ belongs to } C^1(\R^2),
\end{align*}
and satisfies the growth conditions
\begin{align*}
& -a_0^G(t,x)|z|- b_0^G(t,x) \le -G(t,x,z) \le a_1^G |z|^r
+ b_1^G(t,x),\\
& |\partial_z G(t,x,z)|\le a_2^G |z|^{r-1} + b_2^G(t,x),
\end{align*}
for every $(t,x,z)\in[0,T]\times\partial_S\Om\times\R^2$.
Here $r$ is an exponent related to the trace operators on
Sobolev spaces: if $p<2$, then we suppose
that $p \le r \le \frac{p}{2-p}$, while if $p \ge 2$,
we suppose $p \le r$. Moreover
$a_1^G \ge 0$, $a_2^G \ge 0$ are two nonnegative constants,
and
$a_0^G\in L^\infty([0,T]; L^{r'}(\partial_S\Om))$,
$b_0^G, b_1^G \in C^0([0,T];L^1(\partial_S\Om))$,
and $b_2^G\in C^0([0,T];L^{r'}(\partial_S\Om))$ are
nonnegative functions with $r':=r/(r-1)$
\par
We assume that for every $(t,z) \in [0,T] \times \R^2$
\begin{align*}
& G(t,x,z) = G(0,x,z) + \int_0^t \dot{G}(s,x,z) ds
\quad\quad \text{ for }
\hs^1\text{{-}a.e. }x\in\partial_S\Om,\\
& \partial_z G(t,x,z) = \partial_z G(0,x,z) +
\int_0^t \partial_z \dot{G}(s,x,z) ds \quad\quad
\text{ for }\hs^1\text{{-}a.e. }x\in\partial_S\Om,
\end{align*}
where $\dot{G}: [0,T]\times\partial_S\Om\times\R^2 \to \R$
is such that
\begin{align*}
& \text{ for all }z \in \R^2:\, (t,x)\to \dot{G}(t,x,z)
\text{ is } \leb^1\times\hs^1 \text{{-}measurable, } \\
& \text{ for all } (t,x)\in [0,T]\times\partial_S\Om:\,
z \to \dot{G}(t,x,z)
\text{ belongs to } C^1(\R^2),
\end{align*}
and satisfies the the growth conditions
\begin{align*}
& |\dot{G}(t,x,z)|\le a_3^G(t)|z|^r + b_3^G(t,x), \\
& |\partial_z\dot{G}(t,x,z)|\le a_4^G(t)|z|^{r-1} +
b_4^G(t,x)
\end{align*}
for all $(t,x,z)\in[0,T]\times\partial_S\Om\times\R^2$.
Here
$a_3^G,\,a_4^G \in L^1([0,T])$,
$b_3^G\in L^1([0,T];L^1(\partial_S\Om))$ and
$b_4^G \in L^1([0,T];L^{r'}(\partial_S \Om))$ are
nonnegative functions.
\par
By the previous assumptions, the following functionals on
$L^r(\partial_S \Om;\R^2)$
\begin{equation}
\label{tracener}
\gs(t)(u):= \int_{\partial_S\Om} G(t,x,u(x)) \,d\hs^1(x),\quad
\dot{\gs}(t)(u):=\int_{\partial_S\Om} \dot{G}(t,x,u(x)) \,d\hs^1(x)
\end{equation}
are well defined. For every $t \in [0,T]$ we have that
$\gs(t)$ is of class $C^1$ on
$L^r(\partial_S \Om;\R^2)$ and its differential
is given by
\begin{equation*}
\langle\partial\gs(t)(u),v\rangle =
\int_{\partial_S\Om} \partial_z G(t,x,u(x)) v(x) \,d\hs^1(x),
\quad\quad
u,v \in L^{r}(\partial_S\Om;\R^2),
\end{equation*}
where $\langle\cdot,\cdot\rangle$ denotes now the duality
pairing between
$L^{r'}(\partial_S\Om;\R^2)$ and $L^{r}(\partial_S\Om;\R^2)$. 
Moreover, $\dot{\gs}(t)$ is of class $C^1$ on 
$L^{r}(\partial_S \Om;\R^2)$, and its differential
is given by
\begin{equation*}
\langle\partial\dot{\gs}(t)(u),v\rangle =
\int_{\partial_S\Om} \partial_z \dot{G}(t,x,u(x)) v(x)
\,d\hs^1(x)
\end{equation*}
for all $u,v \in L^{r}(\partial_S\Om;\R^2)$. Finally we have
\begin{equation*}
\gs(t)(u) = \gs(0)(u) + \int_0^t \dot{\gs}(s)(u) \,ds,
\end{equation*}
\begin{equation*}
\langle\partial \gs(t)(u),v\rangle =
\langle\partial \gs(0)(u),v\rangle +
\int_0^t \langle\partial \dot{\gs}(s)(u),v\rangle
\,ds,
\end{equation*}
for every $u,v\in L^{r}(\partial_S\Om;\R^2)$.
\par
Let $\Om_S \subseteq \Om\setminus \OmBb$ be open
with Lipschitz boundary, and such
that $\partial_S \Om \subseteq \partial \Om_S$; the trace
operator from $W^{1,p}(\Om_S;\R^2)$
into $L^{r}(\partial \Om_S;\R^2)$ is then compact, and
so there exists a constant
$\gamma_S>0$ such that
\begin{equation}
\label{sf1}
\|u\|_{r,\partial_S \Om} \le
\gamma_S ( \|\nabla u\|_{p, \Om_S} +\|u\|_{p, \Om_S})
\end{equation}
for every $u \in W^{1,p}(\Om_S;\R^2)$.
By the previous assumptions, we have that
there exist six nonnegative constants
$\alpha_0^{\gs}$, $\alpha_1^{\gs}$, $\alpha_2^{\gs}$, $\beta_0^{\gs}$,
$\beta_1^{\gs}$,
$\beta_2^{\gs}$ and four nonnegative functions $\alpha_3^{\gs}$,
$\alpha_4^{\gs}$, $\beta_3^{\gs}$,
$\beta_4^{\gs} \in L^1([0,T])$, such that
\begin{align}
\nonumber
& -\alpha_0^G\|u\|_{r,\partial_S\Om} - \beta_0^G
\le -\gs(t)(u)\le
\alpha_1^{\gs} \|u\|^{r}_{r,\partial_S\Om}+ \beta_1^G,\\
\label{gbou2}
& |\langle\partial \gs(t)(u),v\rangle| \le
(\alpha_2^{\gs}\|u\|^{r-1}_{r,\partial_S\Om} +
\beta_2^{\gs})
\|v\|_{r,\partial_S\Om},\\
\label{gbou3}
& |\dot{\gs}(t)(u)| \le \alpha_3^{\gs}(t)\|u\|^{r}_{r,\partial_S\Om} +
\beta_3^{\gs}(t),\\
\nonumber
& |\langle\partial \dot{\gs}(t)(u),v\rangle| \le
(\alpha_4^{\gs}(t)\|u\|^{r-1}_{r, \partial_S\Om} + \beta_4^{\gs}(t))
\|v\|_{r,\partial_S\Om}
\end{align}
for every $t \in [0,T]$ and $u,v \in L^r(\partial_S \Om;\R^2)$.
\par\medskip\noindent
{\bf Configurations with finite energy.}
The deformations on the boundary $\partial_D \Om$ are given by
(the traces of) functions
$g \in W^{1,p}(\Om;\R^2)\cap L^q(\Om;\R^2)$, where $p,q$ are the
exponents in \eqref{scisandb2} and
\eqref{hpf4} respectively. Given a crack $\Gamma \in \radm$ and 
a boundary deformation $g$,
the set of {\it admissible deformations with finite energy} 
relative to $(g,\Gamma)$
is defined by
\begin{equation*}
AD(g,\Gamma):= \{u \in GSBV^p_q(\Om;\R^2):
S(u) \tsub \Gamma, u=g\,
\;\hs^1\text{{-}a.e. on }\partial_D\Om\setminus\Gamma\},
\end{equation*}
where we recall that
$$
GSBV^p_q(\Om;\R^2):=GSBV^p(\Om;\R^2) \cap L^q(\Om;\R^2),
$$
and the equality $u=g$ on $\partial_D\Om\setminus\Gamma$ is
intended in the sense of traces
(see \cite[Section 2]{DMFT}).
\par
Note that if $u \in GSBV^p_q(\Om;\R^2)$, then
$\ws(u)<+\infty$ and
$|\fs(t)(u)|<+\infty$ for all $t \in [0,T]$. Moreover
since $\Gamma \in \radm$ and $S(u) \tsub \Gamma$,
we have that
$u \in W^{1,p}(\Om_S;\R^2) \cap L^q(\Om_S;\R^2)$ so that
$\gs(t)(u)$ is well defined and
$|\gs(t)(u)|<+\infty$ for all $t \in [0,T]$.
Notice that  there exists always a deformation without
crack which satisfies the boundary condition,
namely the function $g$ itself.
\par\medskip\noindent
{\bf The total energy.}
For every $t\in [0,T]$, the total energy relative to the
configuration $(u,\Gamma)$ with $u \in AD(g,\Gamma)$
is given by
\begin{equation}
\label{totalener}
\Eub(t)(u,\Gamma):=\Eb(t)(u)+\Es(\Gamma),
\end{equation}
where
\begin{equation}
\label{elener}
\Eb(t)(u):=\ws(u)-\fs(t)(u)-\gs(t)(u),
\end{equation}
and $\ws$, $\fs(t)$, $\gs(t)$ and $\Es$ are defined in
\eqref{bws}, \eqref{bodyenergy},
\eqref{tracener} and \eqref{crackener} respectively.
It turns out that
there exist four constants $\alpha^{\Eub}_0>0$,
$\alpha^{\Eub}_1>0$, $\beta^{\Eub}_0 \ge 0$,
$\beta^{\Eub}_1 \ge 0$ such that
\begin{equation}
\label{totenerbelow}
\Eb(t)(u) \ge \alpha^{\Eub}_0( \|\nabla u\|_p^p+
\|u\|_q^q)- \beta^{\Eub}_0,
\end{equation}
\begin{equation*}
\Eb(t)(u) \le \alpha^{\Eub}_1( \|\nabla u\|_p^p+ \|u\|_q^q+
\|u\|_{r, \partial_S \Om}^r)+ \beta^{\Eub}_1,
\end{equation*}
for every $t \in [0,T]$ and $u \in GSBV^p_q(\Om;\R^2)$.
\par\medskip\noindent
{\bf The time dependent boundary deformations.}
We will consider boundary deformations $g(t)$ such that
$$
t \to g(t) \in \,AC([0,T];
W^{1,p}(\Om;\R^2) \cap L^q(\Om;\R^2)),
$$
so that
$$
t \to \dot{g}(t) \in L^1([0,T];
W^{1,p}(\Om;\R^2) \cap L^q(\Om;\R^2)),
$$
and
$$
t \to \nabla\dot{g}(t)\in L^1([0,T];L^p(\Om;\msd)).
$$
\par\medskip\noindent
{\bf The existence result.}
Let $\Gamma_0 \in \radm$ be a preexisting crack.
The next Theorem proved in \cite{DMFT} establishes
the existence of a quasistatic evolution with preexisting crack
$\Gamma_0$.

\begin{theorem}
\label{qse}
Let $\Gamma_0 \in \radm$ be a preexisting crack.
Then there exists
a quasistatic evolution with preexisting
crack $\Gamma_0$ and boundary deformation $g(t)$, i.e.,
there exists a function $t \to (u(t),\Gamma(t))$ from
$[0,T]$ to $GSBV^p_q(\Om;\R^2)\times \radm$
with the following properties:
\begin{itemize}
\item[]
\item[(a)] $(u(0),\Gamma(0))$ is such that
$$
\eub(0)(u(0),\Gamma(0)) = \min \{\eub(0)(v,\Gamma): v \in
AD(g(0),\Gamma), \Gamma_0 \tsub \Gamma\};
$$
\item[]
\item[(b)]
$u(t) \in AD(g(t),\Gamma(t))$  for all $t\in [0,T]$;
\item[]
\item[(c)]
{\it irreversibility}:\,\,
$\Gamma_0 \tsub \Gamma(s)\tilde{\subset}\Gamma(t)$ whenever
$0\le s<t \le T$;
\item[]
\item[(d)]
{\it static equilibrium}:\,\,
for all  $t\in [0,T]$
$$
\eub(t)(u(t),\Gamma(t)) = \min \{\eub(t)(v,\Gamma): v \in
AD(g(t),\Gamma), \Gamma(t) \tsub \Gamma\};
$$
\item[]
\item[(e)]
{\it nondissipativity}:\,\,
the function $t \to E(t):=\eub(t)(u(t),\Gamma(t))$ is
absolutely continuous on  $[0,T]$,
and for a.e.  $t\in [0,T]$
\begin{eqnarray}
\label{nondissqse}
\dot{E}(t) = \langle\partial\ws(\nabla u(t)),
\nabla \dot{g}(t)\rangle
-\langle\partial \fs(t)(u(t)),\dot{g}(t)\rangle-
\dot{\fs}(t)(u(t)) \\
\nonumber
-\langle\partial \gs(t)(u(t)),\dot{g}(t)\rangle-
\dot{\gs}(t)(u(t)).
\end{eqnarray}
\end{itemize}
\end{theorem}

The next theorem gives a compactness and lower
semicontinuity result with respect
to weak convergence in $GSBV^p_q(\Om,\R^2)$ which
will be often used in the next sections.

\begin{theorem}
\label{gsbvcomp}
Let $t_k \in [0,T]$ with $t_k \to t$, and let
$(u_k)\subset GSBV^p_q(\Om;\R^2)$,
$C \in ]0;+\infty[$ such that $S(u_k) \tsub \Omb_B$ and
$$
\Eb(t_k)(u_k)+\Es(S(u_k))  \le C,
$$
where $\Eb$ and $\Es$ are defined as in \eqref{elener} and
\eqref{crackener}.
Then there exists a subsequence $(u_{k_h})_{h \in \N}$ converging to
some $u$ weakly
in $GSBV^p_q(\Om;\R^2)$ such that $S(u) \tsub \Omb_B$,
$$
\Eb(t)(u) \le \liminf_{h\to \infty} \Eb(t_{k_h})(u_{k_h}),
$$
and
$$
\Es(S(u)) \le  \liminf_{h\to \infty} \Es(S(u_{k_h})).
$$
\end{theorem}

\begin{proof}
By \eqref{totenerbelow} and \eqref{ik3}, we have that there
exists $C' \in ]0,+\infty[$ such that
$$
\|\nabla u_k\|_p^p+ \|u_k\|_q^q +\hs^1(S(u_k)) \le C'.
$$
Then we can apply  Theorem \ref{GSBVcompact} with
$g(x,u_k)=|u_k|^q$, obtaining a subsequence $(u_{k_h})_{h \in \N}$
and $u \in GSBV^p(\Om;\R^2)$
such that \eqref{gsbvconv} holds: in particular we may assume
that $u_{k_h} \to u$ pointwise a.e..
We have $u_{k_h} \to u$ strongly in $L^1(\Om;\R^2)$, and by
Fatou's Lemma we have that
$u \in L^q(\Om;\R^2)$ so that $u \in GSBV^p_q(\Om;\R^2)$.
We conclude $u_{k_h} \weak u$ weakly
in $GSBV^p_q(\Om;\R^2)$. By \cite[Theorem 3.7]{A2} we
have that
$$
\Es(S(u)) \le \liminf_h \Es(S(u_{k_h})),
$$
by \cite{K} we have that
$$
\int_\Om W(x,\nabla u)\, dx \le
\liminf_h \int_\Om W(x,\nabla u_{k_h})\, dx,
$$
and by Fatou's Lemma (in the limsup version) we have
$$
\limsup_h \int_\Om F(t_{k_h},x,u_{k_h}(x))\,dx \le
\int_\Om  F(t,x,u(x))\,dx.
$$
Since ${(u_{h_k})}_{|\Om_S}$ is bounded in
$W^{1,p}(\Om_S;\R^2) \cap L^q(\Om_S;\R^2)$,
and the trace operator from $W^{1,p}(\Om_S;\R^2)$
into $L^{r}(\Om_S;\R^2)$ is compact, we get
$$
\lim_h \gs(t_{k_h})(u_{k_h}) = \gs(t)(u),
$$
and so the proof is thus concluded.
\end{proof}

\section{The finite element space and
the transfer of jump}
\label{femspace}
Let $\Om \subseteq \R^2$ be a polygonal set and
let us fix two positive
constants $0<c_1<c_2<+\infty$. By a
{\it regular triangulation} of $\Om$ of
size $\eps$ we intend a finite family of (closed)
triangles $T_i$ such that $\Omb=\bigcup_i T_i$,
$T_i \cap T_j$ is either empty or equal to a common
edge or to a common vertex,
and each $T_i$ contains a ball of diameter
$c_1 \eps$ and is contained in a ball of
diameter $c_2 \eps$.
\par
We indicate by $\rs_\eps(\Om)$ the family of all
regular triangulations of
$\Om$ of size $\eps$.
It turns out that there exist
$0<\vartheta_1 < \vartheta_2 <\pi$
such that for all triangle $T$ belonging to
a regular triangulation
$\tb \in \rs_\eps(\Om)$, the inner angles of $T$ are between
$\vartheta_1$ and $\vartheta_2$. Moreover, every
edge of $T$ has length greater than
$c_1 \eps$ and lower than $c_2 \eps$.
\par
Let $\eps>0$, $\rb_\eps \in \rs_\eps(\Om)$,
and let $a \in ]0,\frac{1}{2}[$.
Let us consider a triangulation $\tb$ nested in
$\rb_\eps$ obtained
dividing each triangle $T \in \rb_\eps$ into four
triangles taking over every edge $[x,y]$ of $T$
a knot $z$ which satisfies
$$
z=tx+(1-t)y, \quad \quad t \in [a,1-a].
$$
We will call these vertices {\it adaptive}, the
triangles obtained gluing these points
{\it adaptive triangles}, and their edges
{\it adaptive edges}.
We denote by $\ts_{\eps,a}(\Om)$ the set of
triangulations $\tb$ constructed in this way.
Note that for all $\tb \in \ts_{\eps,a}(\Om)$
there exist $0<c_1^a<c_2^a<+\infty$ such that
every $T_i \in \tb$ contains a ball of diameter
$c_1^a \eps$ and is contained in a ball
of diameter $c_2^a \eps$. Then there exist
$0<\vartheta_1^a < \vartheta_2^a <\pi$
such that for all $T$ belonging to
$\tb \in \ts_{\eps,a}(\Om)$, the inner angles of $T$ are between
$\vartheta_1^a$ and $\vartheta_2^a$.
Moreover, every edge of $T$ has length greater than
$c_1^a \eps$ and lower than $c_2^a \eps$.
\begin{center}
\label{fig1}
\begin{figure}
\input{fig1.pstex_t}
\caption{}
\end{figure}
\end{center}
\vskip-17pt
\par
From now on for all $\varepsilon>0$ we
fix $\rb_\varepsilon \in \rs_\varepsilon(\Om)$.
We suppose that the brittle part $\Om_B$ and
the region
$\Om_S$ introduced before for the model of
quasistatic growth of fractures
are composed of triangles of
$\rb_\eps$ for all $\eps$.
Moreover we suppose that $\partial_D \Om$ and
$\partial_S \Om$ are composed of edges of $\rb_\eps$ 
for all $\eps$ up to a finite number of points.
Finally, in order to deal with the deformation
at the boundary, it will be useful to consider
$\Om_D$ polygonal such that
$\Om_D \cap \Om=\emptyset$, and
$\partial \Om_D \cap \partial \Om=\partial_D \Om$
up to a finite number of
points. We set
\begin{equation}
\label{omprimo}
\Om':=\Om \cup \Om_D \cup \partial_D \Om,
\end{equation}
and we suppose that the regular triangulation
$\rb_{\eps}$ can be extended to a regular triangulation
of $\Om'$, so that every triangulation $\tb$ in
$\treaom$ can be extended to a triangulation of
$\trea$ considering the middle points of the new edges
as adaptive vertices: we indicate these
extended triangulation with the same symbol
$\tb$.
\par
We consider the following discontinuous finite element space.
We indicate by $\afeaom$ the set of all $u:\Om \to \R^2$
such that there exists a triangulation
$\tb(u) \in \treaom$ nested in $\rb_\varepsilon$ with
$u$ affine on every triangle
$T \in \tb(u)$. For every $u \in \afeaom$,
we indicate by $S(u)$ the family of edges of
$\tb(u)$ inside $\Om$ across which $u$ is discontinuous.
Notice that $u \in SBV(\Om;\R^2)$
and that the notation is consistent with the usual one
employed in the
theory of functions of bounded variation.
Let us set
\begin{equation}
\label{affini}
\afeom:=
\{u:\Om \to \R^2\,:\,u \text{ is continuous and affine
on each triangle }T \in \rb_\eps\}.
\end{equation}
\par
The discretization of the problem will be carried
out using the space
\begin{equation}
\label{fespace}
\afeaomb:=\{u \in \afeaom\,:\,
S(u) \subseteq \OmBb\}.
\end{equation}
Given any $g \in \afeom$, for every
$u \in \afeaomb$ let
\begin{equation}
\label{jump*}
S_D^g(u):= \{x \in \partial_D \Om\,:\,
u(x) \not= g(x) \},
\end{equation}
that is $S_D^g(u)$ denotes the set of edges
of $\partial_D \Om$
at which the boundary condition is not
satisfied. For every $u \in \afeaomb$,
let us also set
\begin{equation}
\label{newjump}
\Sg{g}{u}:=S(u) \cup S_D^g(u).
\end{equation}

An essential tool in the approximation result of this paper
is Proposition \ref{piecetransf2} which
generalizes the piecewise affine transfer of jump
\cite[Proposition 5.1]{GP} to the case of vector valued 
functions with bulk energy $\Eb$
and surface energy $\Es$ of the form \eqref{elener} and
\eqref{crackener} respectively. 
\par
In order to deal with the surface energy $\Es$
we will need the
following geometric construction.
Let $S \subseteq \Om$ be a segment and let us suppose
that $S$ intersects the edges
of $\rb_\eps$ at most in one point for all $\eps>0$.
Let $a \in ]0,\frac{1}{2}[$,
and let $P=S \cap \zeta$, where $\zeta=[x,y]$ is an
edge of $\rb_\eps$: we indicate with $\pi_a(P)$
the projection of $P$ on the segment
$\{tx+(1-t)y\,:\,t \in [a,1-a]\}$.
The {\it interpolating curve} $S_{\eps,a}$ of
$S$ in $\rb_\eps$ with parameter $a$
is given connecting all the $\pi_a(P)$'s belonging to
the same triangle of $\rb_\eps$ (see Figure 2).
\begin{center}
\label{fig2}
\begin{figure}
\input{fig2.pstex_t}
\caption{}
\end{figure}
\end{center}

\begin{lemma}
\label{appes}
Under the previous assumptions, there exists a function
$\eta(a)$ independent of $S$ with
$\eta(a) \to 0$ as $a \to 0$ such that
\begin{equation*}
\limsup_{\eps \to 0} |\Es(S_{\eps,a})-\Es(S)|
\le \eta(a) \Es(S),
\end{equation*}
where $\Es$ is defined in \eqref{crackener}.
\end{lemma}

\begin{proof}
By \eqref{ik3}, we have that there exist $\omega$
and $K_3>0$ such that for all $x_1,x_2 \in \Omb$ and
$|\nu_1|=|\nu_2|=1$
$$
|k(x_1,\nu_1)-k(x_2,\nu_2)| \le \omega(|x_1-x_2|)+
K_3|\nu_1-\nu_2|,
$$
where $\omega\,:\,]0,+\infty[ \to ]0,+\infty[$ is
a decreasing function
such that $\omega(s) \to 0$ as $s \to 0$.
Let $T \in \rb_\eps$ be such that
$T \cap S \not= \emptyset$, and let us choose
$x_T \in T \cap S$ and
$x_T^{\eps,a} \in T \cap S_{\eps,a}$. Let $c_2>0$
denote the characteristic constant of $\rb_\eps$
such that
every $T \in \rb_\eps$ is contained in a ball
of diameter $c_2 \eps$.
Then we have
\begin{multline*}
\left| \int_{S_{\eps,a} \cap T} k(x,\nu_T^{\eps,a})
\,d\hs^1-
\int_{S \cap T} k(x,\nu_T) \,d\hs^1 \right| \\
\le \left|\int_{S_{\eps,a} \cap T}
k(x_T^{\eps,a},\nu_T^{\eps,a}) \,d\hs^1-
\int_{S \cap T} k(x_T,\nu_T) \,d\hs^1 \right|
+\omega(c_2 \eps) \hs^1(S_{\eps,a} \cap T)+
\omega(c_2 \eps) \hs^1(S \cap T) \\
\le
\left|k(x_T^{\eps,a},\nu_T^{\eps,a})
\hs^1(S_{\eps,a} \cap T)- k(x_T,\nu_T)
\hs^1(S \cap T)\right|
+\omega(c_2 \eps) \left[ \hs^1(S_{\eps,a} \cap T)
+\hs^1(S \cap T) \right],
\end{multline*}
where $\nu_T^{\eps,a},\nu_T$ are the (constant) normal
to $S_{\eps,a} \cap T$ and $S \cap T$ respectively.
We have
\begin{multline*}
\left|k(x_T^{\eps,a},\nu_T^{\eps,a})
\hs^1(S_{\eps,a} \cap T)- k(x_T,\nu_T)
\hs^1(S \cap T)\right|
\\
\le k(x_T^{\eps,a},\nu_T^{\eps,a})
\left|\hs^1(S_{\eps,a} \cap T)-
\hs^1(S \cap T) \right|+
\left|k(x_T^{\eps,a},\nu_T^{\eps,a}) -k(x_T,\nu_T)\right|
\hs^1(S \cap T) \\
\le K_2 \left|\hs^1(S_{\eps,a} \cap T)-
\hs^1(S \cap T) \right|+
\omega(|x_T^{\eps,a}-x_T|)\hs^1(S \cap T)
+K_3|\nu_T^{\eps,a}-\nu_T|\hs^1(S \cap T),
\end{multline*}
where $K_2$ is defined in \eqref{ik3}.
We are now ready to conclude: in fact, following
\cite[Lemma 5.2.2]{N}, we can choose the orientation
of $\nu_T^{\eps,a}$
in such a way that
$$
|\nu_T^{\eps,a}-\nu_T|\hs^1(S \cap T) \le  D_2 a \eps,
\quad \quad
\left|\hs^1(S_{\eps,a} \cap T)-\hs^1(S \cap T) \right|
\le D_1 a \eps,
$$
with $D_1,D_2>0$ independent of $T,\eps,a$.
Then, summing up the preceding
inequalities, recalling that the number of triangles of
$\rb_\eps$ intersecting
$S$ is less than $\tilde{c} \eps^{-1} \hs^1(S)$ for
$\eps$ small enough, with $\tilde{c}$
independent of $S$ and $\eps$ (see for example
\cite[Lemma 2.5]{GP}),
we obtain
$$
\limsup_{\eps \to 0} |\Es(S_{\eps,a})-\Es(S)|
\le \rho(a) \hs^1(S),
$$
where $\rho(a):=\tilde{c}(K_2 D_1+ K_3 D_2)a$.
In view of \eqref{ik3}, we conclude that
$$
\limsup_{\eps \to 0} |\Es(S_{\eps,a})-\Es(S)|
\le K_1^{-1} \rho(a) \Es(S),
$$
and so the proof is concluded choosing
$\eta(a):=K_1^{-1} \rho(a)$.
\end{proof}

For all $u \in GSBV^p_q(\Om;\R^2)$ and for all
$g \in W^{1,p}(\Om;\R^2) \cap L^q(\Om;\R^2)$, let us set
\begin{equation}
\label{jumpg}
\Sg{g}{u}:=S(u) \cup \{x \in \partial_D \Om\,:\,
u(x) \not= g(x)\},
\end{equation}
where the inequality is intended in the sense of traces.
We are now in a position to state the piecewise affine transfer
of jump proposition in our setting.

\begin{proposition}
\label{piecetransf2}
Let $a \in ]0,\frac{1}{2}[$, and for all $i=1,\dots,m$
let
$$
u^i_\eps \in \afeaomb,
\quad\quad
u^i \in GSBV^p_q(\Om;\R^2)
$$
be such that
$$
u^i_\eps \weak u^i \quad\quad \text{ weakly in }
GSBV^p_q(\Om;\R^2).
$$
Let moreover $g^i_\eps,h_\eps \in \afeom$,
$g^i,h \in W^{1,p}(\Om;\R^2) \cap L^q(\Om;\R^2)$
be such that
$$
g^i_\eps \to g^i, \quad h_\eps \to h \quad\quad
\text{ strongly in }W^{1,p}(\Om;\R^2) \cap L^q(\Om;\R^2).
$$
Then for every $v \in GSBV^p_q(\Om;\R^2)$ with
$S(v) \tsub \OmBb$,
there exists $v_\eps \in \afeaomb$ such that
\begin{equation*}
\nabla v_\eps \to \nabla v
\quad\quad \mbox{ strongly in }L^p(\Om;\msd),
\end{equation*}
\begin{equation*}
v_\eps \to v
\quad\quad \mbox{ strongly in }L^q(\Om;\R^2),
\end{equation*}
and such that
\begin{equation*}
\limsup_{\eps \to 0} \Es \left (\Sg{h_\eps}{v_\eps}
\setminus \bigcup_{i=1}^m \Sg{g^i_\eps}{u^i_\eps}\right)
\le \mu(a)
\Es \left( \Sg{h}{v}
\setminus \bigcup_{i=1}^m \Sg{g^i}{u^i} \right),
\end{equation*}
where $\mu(a)$ depends only on $a$, $\mu(a) \to 1$ as
$a \to 0$, and $\Es$ is defined in \eqref{crackener}.
In particular for all $t \in [0,T]$ and for all
$t_\eps \to t$ we have
\begin{equation*}
\Eb(t_\eps)(v_\eps) \to \Eb(t)(v),
\end{equation*}
where $\Eb$ is defined in \eqref{elener}.
\end{proposition}

The proof of Proposition \ref{piecetransf2} can be obtained from
that of \cite[Proposition 5.1]{GP} taking into account
the following modifications. We can consider $v$ scalar valued
since vector valued maps can be easily dealt componentwise.
Even if the surface energy is of the form
\eqref{crackener}, we can still restrict ourselves to the case
in which $v$ has piecewise linear jumps outside a suitable
neighborhood of $\bigcup_{i=1}^m \Sg{g^i}{u^i}$ by using the
density result of \cite{CT}. In order to approximate the piecewise
linear jumps, we use Lemma \ref{appes}. Finally the fact that
$p \not= 2$ prevents us from considering the piecewise jumps
as union of disjoint segments: we overcome this difficulty
choosing $v_\eps=0$ in the regular triangles which contain
the intersection points, and then interpolating $v$ outside
as in \cite[Proposition 5.1]{GP}.

\section{Preexisting cracks and their approximation}
\label{secincrack}
In Section \ref{approx}, we will need to approximate the
surface energy of a given preexisting crack $\Gamma^0$.
We take the initial crack in the class
\begin{multline}
\label{gammaom}
\gadm:=\{\Gamma \tsub \OmBb\,:\, \hs^1(\Gamma)<+\infty,\,
\Gamma=\Sg{h}{z}  \\
\mbox{ for some }
h \in W^{1,p}(\Om;\R^2) \cap L^q(\Om;\R^2)
\text{ and }z \in GSBV^p_q(\Om;\R^2)
\}.
\end{multline}
Notice that it is not restrictive to assume $h \equiv 0$.
We take as discretization of $\gadm$ the following class
\begin{equation}
\label{gammaomea}
\gadmea:=\{\Gamma \tsub \OmBb\,:\, \hs^1(\Gamma)<+\infty,\,
\Gamma=\Sg{0}{z} \mbox{ for some }
z \in \afeaomb\}.
\end{equation}
We have the following approximation result.

\begin{proposition}
\label{gammazero}
Let $\Gamma^0 \in \gadm$. Then for every
$\eps>0$ and $a \in ]0,\frac{1}{2}[$
there exists $\Gamma^0_{\eps,a} \in \gadmea$ such that
$$
\lim_{\eps,a \to 0}
\Es(\Gamma^0_{\eps,a})=\Es(\Gamma^0),
$$
where $\Es$ is defined in \eqref{crackener}.
\par
Moreover let $g_{\eps} \in \afeom$,
$g \in W^{1,p}(\Om;\R^2) \cap L^q(\Om;\R^2)$ be such that
as $\eps \to 0$
$$
g_{\eps} \to g \text{ strongly in }
W^{1,p}(\Om;\R^2) \cap L^q(\Om;\R^2),
$$
and let us consider
$$
F_{\eps,a}(v):=
\begin{cases}
\Eb(0)(v)+\Es \left( \Sg{g_{\eps}}{v}
\setminus  \Gamma^0_{\eps,a} \right)
\quad &\mbox{if }v \in \afeaomb,\\
+\infty \quad &\mbox{otherwise in } L^1(\Om;\R^2),
\end{cases}
$$
and
$$
F(v):=
\begin{cases}
\Eb(0)(v)+\Es \left( \Sg{g}{v} \setminus \Gamma^0 \right)
\quad &\mbox{if }v \in GSBV^p_q(\Om;\R^2),
S(v) \tsub \OmBb, \\
+\infty \quad &\mbox{otherwise in } L^1(\Om;\R^2),
\end{cases}
$$
where $\Eb$ is defined in \eqref{elener}.
Then the family $(F_{\eps,a})$ $\Gamma${-}converges to
$F$ in the strong topology of $L^1(\Om;\R^2)$
as $\eps \to 0$ and $a \to 0$.
\end{proposition}

\begin{proof}
Let us consider $\Gamma^0 \in \gadm$ with
$\Gamma^0=\Sg{0}{z}$ for some $z \in GSBV^p_q(\Om;\R^2)$.
Then by Proposition \ref{piecetransf2}
for every $\eps>0$ and $a\in (0,\frac{1}{2})$,
there exists $\tilde{z}_{\eps,a} \in \afeaom$
such that for $\eps \to 0$ and for all $a$
$$
\nabla \tilde{z}_{\eps,a} \to \nabla z
\quad\quad
\text{ strongly in }L^p(\Om;\msd),
$$
$$
\tilde{z}_{\eps,a} \to z
\quad\quad
\text{ strongly in }L^q(\Om;\R^2),
$$
and
$$
\limsup_{\eps \to 0}\Es(\Sg{0}{\tilde{z}_{\eps,a}}) \le \mu(a)\Es(\Sg{0}{z})
$$
with $\mu(a) \to 1$ as $a \to 0$,
where $\Es$ is defined in \eqref{crackener}.
Let $a_i \searrow 0$, and let $\eps_i \searrow 0$ be
such that for all $\eps \le \eps_i$
$$
\Es(\Sg{0}{\tilde{z}_{\eps,a_i}}) \le \mu(a_i)\Es(\Sg{0}{z})+a_i,
$$
and
$$
\|\nabla \tilde{z}_{\eps,a_i} -\nabla z\|_{L^p(\Om;\msd)} \le a_i,
\quad\quad
\|\tilde{z}_{\eps,a_i} -z\|_{L^q(\Om;\R^2)} \le a_i.
$$
Setting
\begin{equation*}
z_{\eps,a}:=
\begin{cases}
\tilde{z}_{\eps,a_i}
& \eps_{i+1}<\eps \le \eps_i,\,\,a \le a_i, \\
\tilde{z}_{\eps,a_{j-1}}
& \eps_{i+1}<\eps \le \eps_i,\,\, a_j< a \le a_{j-1},\,j \le i,
\end{cases}
\end{equation*}
we have that
$$
\lim_{\eps,a \to 0}\nabla z_{\eps,a}=\nabla z
\quad\quad\text{ strongly in }L^p(\Om;\msd),
$$
$$
\lim_{\eps,a \to 0}z_{\eps,a}=z
\quad\quad\text{ strongly in }L^q(\Om;\R^2),
$$
and
$$
\limsup_{\eps,a \to 0}\Es(\Sg{0}{z_{\eps,a}})
\le \Es(\Sg{0}{z}).
$$
Since by Theorem \ref{gsbvcomp} we have $\Es(\Sg{0}{z_{\eps,a}})
\le \liminf_{\eps,a \to 0}\Es(\Sg{0}{z_{\eps,a}})$, we conclude that
$$
\lim_{\eps,a \to 0}\Es(\Sg{0}{z_{\eps,a}})
=\Es(\Sg{0}{z}).
$$
Let us set for every $\eps,a$
\begin{equation*}
\Gamma^0_{\eps,a}:=\Sg{0}{z_{\eps,a}}.
\end{equation*}
We have that
\begin{equation*}
\lim_{\eps,a \to 0}\Es(\Gamma^0_{\eps,a})=
\Es(\Gamma^0).
\end{equation*}
\par
Let us come to the second part of the proof.
Let us consider $(\eps_n,a_n)_{n \in \N}$ such that
$\eps_n \to 0$ and $a_n \to 0$. If we prove that
$(F_{\eps_n,a_n})_{n \in \N}$
$\Gamma${-}converges to $F$
in the strong topology of $L^1(\Om;\R^2)$, the proposition
is proved since the sequence is arbitrary. Since we can
reason up to subsequences, it is not restrictive to assume
$a_n \searrow 0$.
\par
Let us start with the $\Gamma${-}limsup inequality
considering $v \in GSBV^p_q(\Om;\R^2),$ with
$S(v) \subseteq \OmBb$. For any $n$ fixed,
by Proposition \ref{piecetransf2} there
exists $\tilde{v}_{\eps,a_n} \in \afeanom$ 
such that for $\eps \to 0$
$$
\nabla \tilde{v}_{\eps,a_n} \to \nabla v
\quad\quad
\text{ strongly in } L^p(\Om;\msd),
$$
$$
\tilde{v}_{\eps,a_n} \to v
\quad\quad
\text{ strongly in } L^q(\Om;\R^2),
$$
and such that
$$
\limsup_{\eps \to 0}
\Es(\Sg{g_\eps}{\tilde{v}_{\eps,a_n}} \setminus
\Gamma^0_{\eps,a_n}) \le
\mu(a_n) \Es(\Sg{g_\eps}{v} \setminus \Gamma^0)
$$
with $\mu(a) \to 1$ as $a \to 0$.
For every $m \in \N$ let $\eps^m$ be
such that for all $\eps \le \eps^m$
$$
\Es(\Sg{g_\eps}{\tilde{v}_{\eps,a_m}} \setminus
\Gamma^0_{\eps,a_m}) \le
\mu(a_m)\Es(\Sg{g}{v} \setminus \Gamma^0)+a_m,
$$
and
$$
\|\nabla \tilde{v}_{\eps,a_m} -\nabla v\|_{L^p(\Om;\msd)} \le a_m,
\quad\quad
\|\tilde{v}_{\eps,a_m} -v\|_{L^q(\Om;\R^2)} \le a_m.
$$
We can assume $\eps^m \searrow 0$. Setting
\begin{equation*}
v_{\eps_n,a_n}:=
\begin{cases}
\tilde{v}_{\eps_n,a_m}
& \eps^{m+1}<\eps_n \le \eps^m,\,\, n \ge m, \\
\tilde{v}_{\eps_n,a_n}
& \eps^{m+1}<\eps_n \le \eps^m,\,\, n<m,
\end{cases}
\end{equation*}
we have that
$$
\lim_n\nabla v_{\eps_n,a_n}=\nabla v
\quad\text{ strongly in }L^p(\Om;\msd),
$$
$$
\lim_n v_{\eps_n,a_n}=v
\quad\text{ strongly in }L^q(\Om;\R^2),
$$
and
$$
\limsup_n\Es(\Sg{g_\eps}{v_{\eps_n,a_n}}
\setminus \Gamma^0_{\eps,a}) \le
\Es(\Sg{g}{v} \setminus \Gamma^0).
$$
Then we get
\begin{align*}
\limsup_n
F_{\eps_n,a_n}(v_{\eps_n,a_n}) &\le
\limsup_n \Eb(0)(v_{\eps_n,a_n})+
\limsup_n
\Es(\Sg{g_{\eps_n}}{v_{\eps_n,a_n}} \setminus \Gamma^0_{\eps_n,a_n}) \\
&\le \Eb(0)(v)+\Es(\Sg{g}{v} \setminus \Gamma^0)=F(v),
\end{align*}
so that the $\Gamma${-}limsup inequality holds.
\par
Let us come to the $\Gamma${-}liminf inequality.
Let $v_n,v \in L^1(\Om;\R^2)$ be
such that $v_n \to v$ strongly in $L^1(\Om;\R^2)$ and
$\liminf_n F_{\eps_n,a_n}(v_n)<+\infty$.
By Theorem \ref{gsbvcomp}, we have
$v \in GSBV^p_q(\Om;\R^2)$ with
$S(v) \tsub \OmBb$ and
$$
\Eb(0)(v) \le \liminf_n \Eb(0)(v_n).
$$
Let us consider $\Om'$ defined in \eqref{omprimo}.
Let us extend $g_{\eps_n}$ and $g$ to
$W^{1,p}(\Om';\R^2) \cap L^q(\Om';\R^2)$
in such a way that
$g_{\eps_n} \to g$ strongly in
$W^{1,p}(\Om';\R^2) \cap L^q(\Om';\R^2)$,
and let us also extend $v_n,v$ to
$\Om'$ setting $v_n=g_{\eps_n}$ and
$v=g$ on $\Om_D$. We indicate these
extensions with $w_n$ and $w$ respectively.
Let us also set $z_{\eps_n,a_n}=z=0$ on $\Om_D$,
where $z_{\eps_n,a_n}$ and $z$ are such that
$\Gamma^0_{\eps,a}=z_{\eps_n,a_n}$ and $\Gamma^0=S(z)$.
We indicate these extension by $\zeta_{\eps_n,a_n}$ and
$\zeta$ respectively.
Then for every $\eta>0$ we have by Theorem
\ref{gsbvcomp}
$$
\Es(S(w+\eta \zeta)) \le \liminf_n
\Es( S(w_n+\eta \zeta_{\eps_n,a_n})).
$$
Since for a.e. $\eta>0$ we have
$S(w+\eta \zeta)=S(w) \cup S(\zeta)$ and
$S(w_n+\eta \zeta_{\eps_n,a_n})=S(w_n) \cup
S(\zeta_{\eps_n,a_n})$, we deduce that
$$
\Es(\Sg{g}{v} \cup \Gamma^0) \le
\liminf_n \Es(\Sg{g_{\eps_n}}{v_n} \cup
\Gamma^0_{\eps_n,a_n}).
$$
Since by assumption
$\Es(\Gamma^0_{\eps_n,a_n}) \to \Es(\Gamma^0)$,
we conclude that
$$
\Es(\Sg{g}{v} \setminus \Gamma^0) \le
\liminf_n \Es(\Sg{g_{\eps_n}}{v_n} \setminus
\Gamma^0_{\eps_n,a_n}).
$$
We deduce that
$$
\Eb(0)(v)+ \Es(\Sg{g}{v} \setminus \Gamma^0) \le
\liminf_n
\left[ \Eb(0)(v_n)+ \Es(\Sg{g_{\eps_n}}{v_n}
\setminus \Gamma^0_{\eps_n,a_n}) \right]
$$
that is
$$
F(v) \le \liminf_n F_{\eps_n,a_n}(v_n).
$$
The $\Gamma${-}liminf inequality holds, and so
the proof is concluded.
\end{proof}

\section{The discontinuous finite element approximation}
\label{devol}
In this section we construct a discrete approximation
of the quasistatic evolution
of brittle fractures proposed in \cite{DMFT} and
described in the Preliminaries:
the discretization is done both in space and time.
Let us consider
$$
g_\eps \in W^{1,1}([0,T];
W^{1,p}(\Om;\R^2) \cap L^q(\Om;\R^2)),
\quad
g_\eps(t) \in \afeom \text{ for all }t \in [0,T],
$$
where $\afeom$ is defined in \eqref{affini}.
Let $\delta>0$, and
let $N_\delta$ be the largest integer such that
$\delta (N_\delta-1) < T$; we set
$t_i^\delta:=i\delta$ for
$0\le i \le N_\delta-1$,
$t^\delta_{N_\delta}:=T$ and
$g_\eps^{\delta,i}:=g_\eps(t_i^\delta)$.
Let $\Gamma^0 \in \gadmea$ be a preexisting crack in $\Om$,
where $\gadmea$ is defined in \eqref{gammaomea}.

\begin{proposition}
\label{discrevol}
Let $\eps>0$, $a \in ]0,\frac{1}{2}[$ and
$\delta>0$ be fixed. Then for all
$i=0, \ldots, N_\delta$ there exists
$u^{\delta,i}_{\eps,a} \in \afeaomb$
such that, setting
\begin{equation*}
\Gamma^{\delta,i}_{\eps,a}:=
\Gamma^0 \cup \bigcup_{r=0}^i
\Sg{g_\eps^{\delta,r}}{u^{\delta,r}_{\eps,a}},
\end{equation*}
we have for all
$v \in \afeaomb$
\begin{equation}
\label{MSdiscr}
\Eb(0)(u^{\delta,0}_{\eps,a})+
\Es \left(
\Sg{g_\eps^{\delta,0}}{u^{\delta,0}_{\eps,a}}
\setminus \Gamma^0 \right)
\le \Eb(0)(v)+
\Es \left( \Sg{g_\eps^{\delta,0}}{v}
\setminus \Gamma^0 \right),
\end{equation}
and for $1 \le i \le N_\delta$
\begin{equation}
\label{piecemin}
\Eb(t^\delta_i)(u^{\delta,i}_{\eps,a})+
\Es\left( \Sg{g_\eps^{\delta,i}}{u^{\delta,i}_{\eps,a}}
\setminus \Gamma^{\delta, i-1}_{\eps,a} \right)
\le \Eb(t^\delta_{i})(v)+
\Es\left( \Sg{g_\eps^{\delta,i}}{v}
\setminus \Gamma^{\delta, i-1}_{\eps,a} \right).
\end{equation}
\end{proposition}

\begin{proof}
Let $u^{\delta,0}_{\eps,a}$ be a minimum of the
following problem
\begin{equation}
\label{step0}
\min_{u \in \afeaomb}
\left\{ \Eb(0)(u)+
\Es(\Sg{g_\eps^{\delta,0}}{u} \setminus \Gamma^0)\right\}.
\end{equation}
We set $\Gamma^{\delta,0}_{\eps,a}:= \Gamma^0 \cup
\Sg{g_\eps^{\delta,0}}{u^{\delta,0}_{\eps,a}}$.
Recursively, supposing to have constructed
$u^{\delta,i-1}_{\eps,a}$ and
$\Gamma^{\delta,i-1}_{\eps,a}$, let
$u^{\delta,i}_{\eps,a}$ be a minimum for
\begin{equation}
\label{stepj}
\min_{u \in \afeaomb}
\left\{ \Eb(t^\delta_i)(u)+
\Es(\Sg{g_\eps^{\delta,i}}{u}
\setminus \Gamma^{\delta,i-1}_{\eps,a})
\right\}.
\end{equation}
We set $\Gamma^{\delta,i}_{\eps,a}:=
\Gamma^{\delta,i-1}_{\eps,a} \cup
\Sg{g_\eps^{\delta,i}}{u^{\delta,i}_{\eps,a}}$.
It is clear by construction that
\eqref{MSdiscr} and \eqref{piecemin} hold.
\par
Let us prove that problem \eqref{stepj} admits a
solution, problem \eqref{step0} being similar.
Let $(u_n)_{n \in \N}$ 
be a minimizing sequence for problem
\eqref{stepj}: since
$g_\eps^{\delta,i}$ is an admissible test function,
we deduce that for $n$ large
\begin{equation*}
\Eb(t^\delta_i)(u_n)+ \Es \left( \Sg{g_\eps^{\delta,i}}{u_n}
\setminus \Gamma^{\delta,i-1}_{\eps,a} \right) \le
\Eb(t^\delta_i)(g_\eps^{\delta,i})+1.
\end{equation*}
By the lower estimate on the elastic energy \eqref{totenerbelow},
we deduce that for $n$ large
\begin{equation}
\label{boundstepj}
\alpha^{\Eub}_0 \left( \|\nabla u_n\|^p_p
+\|u_n\|^q_q \right)
+\Es \left( \Sg{g_\eps^{\delta,i}}{u_n}
\setminus \Gamma^{\delta,i-1}_{\eps,a} \right)
\le \Eb(t^\delta_i)(g_\eps^{\delta,i})+1+
\beta^\eub_0.
\end{equation}
Let us indicate by $T_n^1, \ldots, T_n^k$ the
triangles of $\tb(u_n)$.
Up to a subsequence, there exists
$\tb=\{T^1, \ldots, T^k\} \in \treaom$ such that for all
$i=1, \ldots, k$ we have $T_n^i \to T^i$ in the Hausdorff metric
(see Section \ref{notprel} for a precise definition). Let us consider
$T^i \in \tb$, and let $\tilde{T}^i$ be contained in
the interior of $T^i$. For $n$
large enough, $\tilde{T}^i$ is contained in the
interior of $T^i_n$; moreover
$(u_n)_{|\tilde{T}^i}$ is affine and by
\eqref{boundstepj} we have
$\int_{\tilde{T}^i} |\nabla u_n|^p \,dx +
\|u_n\|_{L^\infty(\tilde{T}^i;\R^2)}
\le C$, with $C$ independent of $n$.
We deduce that there exists a function $u^i$
affine on $\tilde{T}^i$ such that up to a
subsequence $u_n \to u$ uniformly on $\tilde{T}^i$.
Since $\tilde{T}^i$ is arbitrary,
it turns out that $u^i$ is actually defined on
$T^i$ and
\begin{equation*}
\Eb(t^\delta_i)_{|T^i}(u^i) \le
\liminf_n \Eb(t^\delta_i)_{|T^i}(u_n^i),
\end{equation*}
where $\Eb(t^\delta_i)_{|T^i}$ denotes the restriction
of $\Eb(t^\delta_i)$ to the
maps defined on $T^i$.
Let $u \in \afeaom$ be such that $u=u^i$ on $T^i$ for every
$i=1,\dots,k$: we have
\begin{equation*}
\Eb(t^\delta_i)(u) \le \liminf_n \Eb(t^\delta_i)(u_n).
\end{equation*}
On the other hand, it is easy to see that
$\Sg{g_\eps^{\delta,i}}{u}$ is contained in the Hausdorff limit
of $\Sg{g_\eps^{\delta,i}}{u_n}$, and so $u \in \afeaomb$;
moreover we deduce
\begin{equation*}
\Es \left( \Sg{g_\eps^{\delta,i}}{u} \setminus
\Gamma^{\delta,i-1}_{\eps,a}
\right)
\le
\liminf_n \Es \left( \Sg{g_\eps^{\delta,i}}{u_n} \setminus
\Gamma^{\delta,i-1}_{\eps,a}
\right).
\end{equation*}
We conclude that $u$ is a minimum point for the problem
\eqref{stepj}, so that the proof is concluded.
\end{proof}

The following estimate on the total energy is essential
in order to study the asymptotic behavior of
the discrete evolution as $\delta \to 0$, $\eps \to 0$
and $a \to 0$.
Let us set $u^\delta_{\eps,a}(t):=u^{\delta,i}_{\eps,a}$
for all $t^\delta_i \le t <t^\delta_{i+1}$ and
$i=0,\dots,N_\delta-1$, $u^\delta_{\eps,a}(T)=
u^{\delta,N_\delta}_{\eps,a}$.

\begin{proposition}
\label{dener}
For all $0 \le j \le i \le N_{\delta}$ we have
\begin{multline}
\label{discrenergy}
\eub(t^\delta_i)(u^{\delta,i}_{\eps,a},\Gamma^{\delta,i}_{\eps,a})
 \le
\eub(t^\delta_j)(u^{\delta,j}_{\eps,a},\Gamma^{\delta,j}_{\eps,a})
+\int_{t^\delta_j}^{t^\delta_i}
\langle \partial \ws(\nabla u^\delta_{\eps,a}(\tau)),
\nabla \dot{g}_\eps(\tau) \rangle \,d\tau \\
-\int_{t^\delta_j}^{t^\delta_i}
\dot{\fs}(\tau)(u^\delta_{\eps,a}(\tau))\,d\tau
-\int_{t^\delta_j}^{t^\delta_i}
\langle \partial \fs(\tau)(u^\delta_{\eps,a}(\tau)),
\dot{g}_\eps(\tau) \rangle \,d\tau \\
-\int_{t^\delta_j}^{t^\delta_i}
\dot{\gs}(\tau)(u^\delta_{\eps,a}(\tau))\,d\tau
-\int_{t^\delta_j}^{t^\delta_i}
\langle \partial \gs(\tau)(u^\delta_{\eps,a}(\tau)),
\dot{g}_\eps(\tau) \rangle \,d\tau
+e^\delta_{\eps,a},
\end{multline}
where $e^\delta_{\eps,a} \to 0$ as $\delta \to 0$ uniformly in
$\eps$ and $a$.
\end{proposition}

\begin{proof}
By the minimality property \eqref{piecemin}, comparing
$u^{\delta,i}_{\eps,a}$ with
$u^{\delta,i-1}_{\eps,a}-g_\eps^{\delta,i-1}+g_\eps^{\delta,i}$
we get
\begin{multline}
\label{comparison}
\ws(\nabla u^{\delta,i}_{\eps,a})-
\fs(t^\delta_{i})(u^{\delta,i}_{\eps,a})
-\gs(t^\delta_{i})(u^{\delta,i}_{\eps,a})+
\Es(\Sg{g^\delta_i}{u^{\delta,i}_{\eps,a}}
\setminus \Gamma^\delta_{i-1}) \\
\le \ws(\nabla u^{\delta,i-1}_{\eps,a}-
\nabla  g_\eps^{\delta,i-1}+\nabla  g_\eps^{\delta,i})
-\fs(t^\delta_{i})(u^{\delta,i-1}_{\eps,a}-
g_\eps^{\delta,i-1}+ g_\eps^{\delta,i}) \\
-\gs(t^\delta_{i})(u^{\delta,i-1}_{\eps,a}-
g_\eps^{\delta,i-1}+ g_\eps^{\delta,i}).
\end{multline}
We have
\begin{multline}
\label{ws}
\ws(\nabla u^{\delta,i-1}_{\eps,a}-
\nabla  g_\eps^{\delta,i-1}+\nabla g_\eps^{\delta,i})
= \ws(\nabla u^{\delta,i-1}_{\eps,a}) \\
+\langle \partial
\ws(\nabla u^{\delta,i-1}_{\eps,a}+ \vartheta_{\eps,a}^{\delta,i-1}
(\nabla  g_\eps^{\delta,i}-\nabla  g_\eps^{\delta,i-1})),
\nabla  g_\eps^{\delta,i}-\nabla  g_\eps^{\delta,i-1} \rangle \\
=\ws(\nabla u^{\delta,i-1}_{\eps,a})+
\int_{t^\delta_{i-1}}^{t^\delta_i}
\langle \partial \ws(\nabla u^{\delta}_{\eps,a}(\tau)+
v_{\eps,a}^\delta(\tau)),
\nabla  \dot{g}_\eps(\tau) \rangle \,d\tau,
\end{multline}
where $\vartheta_{\eps,a}^{\delta,i-1} \in ]0,1[$ and
$v_{\eps,a}^\delta(\tau):= \vartheta_{\eps,a}^{\delta,i-1}
(\nabla g_\eps^{\delta,i}-  \nabla g_\eps^{\delta,i-1})$
for all $\tau \in [t^\delta_{i-1}, t^\delta_i[$.
\par
Similarly we obtain
\begin{equation}
\label{fs}
\fs(t^\delta_i)(u^{\delta,i-1}_{\eps,a}
-g_\eps^{\delta,i-1}+  g_\eps^{\delta,i})
= \fs(t^\delta_i)(u^{\delta,i-1}_{\eps,a})+
\int_{t^\delta_{i-1}}^{t^\delta_i}
\langle \partial
\fs(t^\delta_i)(u^{\delta}_{\eps,a}(\tau)+ w_{\eps,a}^\delta(\tau)),
\dot{g}_\eps(\tau) \rangle \,d\tau,
\end{equation}
and
\begin{equation}
\label{gs}
\gs(t^\delta_i)(u^{\delta,i-1}_{\eps,a}
- g_\eps^{\delta,i-1}+ g_\eps^{\delta,i})
= \gs(t^\delta_i)(u^{\delta,i-1}_{\eps,a})+
\int_{t^\delta_{i-1}}^{t^\delta_i}
\langle \partial
\gs(t^\delta_i)(u^{\delta}_{\eps,a}(\tau)+ z_{\eps,a}^\delta(\tau)),
\dot{g}_\eps(\tau) \rangle \,d\tau,
\end{equation}
where $w_{\eps,a}^\delta(\tau):= \lambda_{\eps,a}^{\delta,i-1}
(g_\eps^{\delta,i}-g_\eps^{\delta,i-1})$,
$z_{\eps,a}^\delta(\tau):=\nu_{\eps,a}^{\delta,i-1}
(g_\eps^{\delta,i}-g_\eps^{\delta,i-1})$
for all $\tau \in [t^\delta_{i-1}, t^\delta_i[$,
and $\lambda_{\eps,a}^{\delta,i-1}, \nu_{\eps,a}^{\delta,i-1}
\in ]0,1[$.
\par
Since by \eqref{ac2bis} we have for
$\tau \in [t^\delta_{i-1},t^\delta_i[$
\begin{multline*}
\langle \partial
\fs(t^\delta_i)(u^{\delta}_{\eps,a}(\tau)+ w_{\eps,a}^\delta(\tau)),
 \dot{g}_\eps(\tau) \rangle -
\langle \partial
\fs(\tau)(u^{\delta}_{\eps,a}(\tau)+ w_{\eps,a}^\delta(\tau)),
 \dot{g}_\eps(\tau) \rangle \\
= \int_{\tau}^{t^\delta_i}
\langle \partial
\dot{\fs}(s)(u^{\delta}_{\eps,a}(\tau)+ w_{\eps,a}^\delta(\tau)),
 \dot{g}_\eps(\tau) \rangle \,ds
\end{multline*}
we get by \eqref{fs4}
\begin{multline}
\label{fstime2}
\left|
\langle \partial
\fs(t^\delta_i)(u^{\delta}_{\eps,a}(\tau)+ w_{\eps,a}^\delta(\tau)),
 \dot{g}_\eps(\tau) \rangle -
\langle \partial \fs(\tau)(u^{\delta}_{\eps,a}(\tau)+
w_{\eps,a}^\delta(\tau)),
 \dot{g}_\eps(\tau) \rangle
\right|  \\
\le \int_{\tau}^{t^\delta_i}
| \langle \partial
\dot{\fs}(s)(u^{\delta}_{\eps,a}(\tau)+
w_{\eps,a}^\delta(\tau)), \dot{g}_\eps(\tau) \rangle |\,ds \\
\le \int_{\tau}^{t^\delta_i}
\left[ \alpha_4^\fs(s) \|u^\delta_{\eps,a}(\tau)+
w_{\eps,a}^\delta(\tau)\|^{\dot{q}-1}_{\dot{q}} +\beta_4^\fs(s) \right]
\|\dot{g}_\eps(\tau)\|_{\dot{q}}\,ds
\le \gamma^{\delta,\eps,a}_\fs \|\dot{g}_\eps(\tau)\|_{\dot{q}},
\end{multline}
where
\begin{equation*}
\gamma^{\delta,\eps,a}_\fs:= \max_{1 \le i \le N_\delta}
\left(
\|u^{\delta,i-1}_{\eps,a}+
\lambda_{\eps,a}^{\delta,i-1}
(g_\eps^{\delta,i}-g_\eps^{\delta,i-1})\|_{\dot{q}-1}^{\dot{q}}
\int_{t^\delta_{i-1}}^{t^\delta_i} \alpha_4^\fs(s) \,ds
+\int_{t^\delta_{i-1}}^{t^\delta_i} \beta_4^\fs(s) \,ds
\right).
\end{equation*}
Similarly we obtain
\begin{equation}
\label{gstime}
\left|
\langle \partial \gs(t^\delta_i)(u^{\delta}_{\eps,a}(\tau)+
z_{\eps,a}^\delta(\tau)),
\dot{g}_\eps(\tau) \rangle -
\langle \partial \gs(\tau)(u^{\delta}_{\eps,a}(\tau)+
z_{\eps,a}^\delta(\tau)),
\dot{g}_\eps(\tau) \rangle
\right|
\le \gamma^{\delta,\eps,a}_\gs \| \dot{g}_\eps(\tau)\|_{r, \partial_S \Om},
\end{equation}
where
\begin{equation*}
\gamma^{\delta,\eps,a}_\gs:= \max_{1 \le i \le N_\delta}
\left(
\|u^{\delta,i-1}_{\eps,a}+ \nu_{\eps,a}^{\delta,i-1}
(g_\eps^{\delta,i}-g_\eps^{\delta,i-1})\|^{r-1}_{r, \partial_S \Om}
\int_{t^\delta_{i-1}}^{t^\delta_i} a_4^\gs(s) \,ds
+\int_{t^\delta_{i-1}}^{t^\delta_i} b_4^\gs(s) \,ds
\right).
\end{equation*}
From \eqref{comparison}, taking into account \eqref{ws},
\eqref{fs}, \eqref{gs}, \eqref{fstime2}, \eqref{gstime},
we have
\begin{multline}
\label{denergy}
\eub(t^\delta_i)(u^{\delta,i}_{\eps,a},
\Gamma^{\delta,i}_{\eps,a}) \le
\eub(t^\delta_{i-1})(u^{\delta,{i-1}}_{\eps,a},
\Gamma^{\delta,i-1}_{\eps,a})
+\int_{t^\delta_{i-1}}^{t^\delta_i}
\langle \partial \ws(\nabla u^\delta_{\eps,a}(\tau)
+v_{\eps,a}^\delta(\tau)),
\nabla  \dot{g}_\eps(\tau) \rangle \,d\tau \\
-\int_{t^\delta_{i-1}}^{t^\delta_i}
\dot{\fs}(\tau)(u^\delta_{\eps,a}(\tau)) \,d\tau-
\int_{t^\delta_{i-1}}^{t^\delta_i}
\langle \partial \fs(\tau)(u^\delta_{\eps,a}(\tau)
+w_{\eps,a}^\delta(\tau)),  \dot{g}_\eps(\tau) \rangle \,d\tau \\
-\int_{t^\delta_{i-1}}^{t^\delta_i}
\dot{\gs}(\tau)(u^\delta_{\eps,a}(\tau)) \,d\tau-
\int_{t^\delta_{i-1}}^{t^\delta_i}
\langle \partial
\gs(\tau)(u^\delta_{\eps,a}(\tau)+z_{\eps,a}^\delta(\tau)),
\dot{g}_\eps(\tau) \rangle \,d\tau \\
+\gamma^{\delta,\eps,a}_\fs\int_{t^\delta_{i-1}}^{t^\delta_i}
\| \dot{g}_\eps(\tau)\|_{\dot{q}} \,d\tau
+\gamma^{\delta,\eps,a}_\gs \int_{t^\delta_{i-1}}^{t^\delta_i}
\| \dot{g}_\eps(\tau)\|_{r, \partial_S \Om} \,d\tau.
\end{multline}
Taking now $0 \le j \le i \le N_\delta$, summing in
\eqref{denergy} form $t^\delta_j$ to $t^\delta_i$, we obtain
\begin{multline}
\label{denergy2}
\eub(t^\delta_i)(u^{\delta,i}_{\eps,a},
\Gamma^{\delta,i}_{\eps,a}) \le
\eub(t^\delta_j)(u^{\delta,j}_{\eps,a},
\Gamma^{\delta,j}_{\eps,a})
+\int_{t^\delta_j}^{t^\delta_i}
\langle \partial
\ws(\nabla u^\delta_{\eps,a}(\tau)+v_{\eps,a}^\delta(\tau)),
\nabla\dot{g}_\eps(\tau) \rangle \,d\tau \\
-\int_{t^\delta_j}^{t^\delta_i}
\dot{\fs}(\tau)(u^\delta_{\eps,a}(\tau)) \,d\tau-
\int_{t^\delta_j}^{t^\delta_i}
\langle \partial
\fs(\tau)(u^\delta_{\eps,a}(\tau)+w_{\eps,a}^\delta(\tau)),
\dot{g}_\eps(\tau) \rangle \,d\tau \\
-\int_{t^\delta_j}^{t^\delta_i}
\dot{\gs}(\tau)(u^\delta_{\eps,a}(\tau)) \,d\tau-
\int_{t^\delta_j}^{t^\delta_i}
\langle \partial
\gs(\tau)(u^\delta_{\eps,a}(\tau)+z_{\eps,a}^\delta(\tau)),
\dot{g}_\eps(\tau) \rangle \,d\tau \\
+\gamma^{\delta,\eps,a}_\fs\int_{t^\delta_j}^{t^\delta_i}
\| \dot{g}_\eps(\tau)\|_{\dot{q}} \,d\tau
+\gamma^{\delta,\eps,a}_\gs \int_{t^\delta_j}^{t^\delta_i}
\|\dot{g}_\eps(\tau)\|_{r, \partial_S \Om} \,d\tau.
\end{multline}
Setting
\begin{multline}
\label{opiccolo}
e^\delta_{\eps,a}:=
\int_0^1
|\langle \partial
\ws(\nabla u^\delta_{\eps,a}(\tau)+v_{\eps,a}^\delta(\tau)),
\nabla\dot{g}_\eps(\tau) \rangle
-\langle \partial \ws(\nabla u^\delta_{\eps,a}(\tau)),
\nabla\dot{g}_\eps(\tau) \rangle| \,d\tau \\
+\int_0^1
|\langle \partial
\fs(\tau)(u^\delta_{\eps,a}(\tau)+w_{\eps,a}^\delta(\tau)),
\dot{g}_\eps(\tau) \rangle
-\langle \partial \fs(\tau)(u^\delta_{\eps,a}(\tau)),
\dot{g}_\eps(\tau) \rangle| \,d\tau \\
+\int_0^1
|\langle \partial
\gs(\tau)(u^\delta_{\eps,a}(\tau)+z_{\eps,a}^\delta(\tau)),
\dot{g}_\eps(\tau) \rangle
-\langle \partial \gs(\tau)(u^\delta_{\eps,a}(\tau)),
\dot{g}_\eps(\tau) \rangle| \,d\tau \\
+\gamma^{\delta,\eps,a}_\fs\int_0^1
 \|\dot{g}_\eps(\tau)\|_{\dot{q}} \,d\tau
+\gamma^{\delta,\eps,a}_\gs \int_0^1
 \| \dot{g}_\eps(\tau)\|_{r, \partial_S \Om} \,d\tau,
\end{multline}
from \eqref{denergy2} we formally obtain \eqref{discrenergy}.
Let us prove that $e^\delta_{|eps,a} \to 0$
as $\delta \to 0$ uniformly in $\eps$ and $a$. 
By \eqref{piecemin}, comparing
$u^{\delta,i}_{\eps,a}$ with $g_\eps^{\delta,i}$,
and taking into account \eqref{totenerbelow},
we get for all $i=1,\dots,N_\delta$,
$$
\|\nabla u^{\delta,i}_{\eps,a}\|_p+
\|u^{\delta,i}_{\eps,a}\|_{q}
\le C',
$$
where
$$
C':=\frac{1}{\alpha^\Eub_0} \max_{i=0,\dots,N_\delta}
\left( \Eb(t^\delta_i)(g_\eps^{\delta,i})+ \beta^\Eub_0 \right).
$$
Since $\Om_S$ is Lipschitz, there exists $K_S>0$
depending only on $p,q$ such that
$$
\|u\|_{p,\Om_S} \le K_S(\|\nabla u\|_{p,\Om_S}+\|u\|_{q,\Om_S})
$$
for all $u \in W^{1,p}(\Om_S;\R^2) \cap L^q(\Om_S;\R^2)$.
Taking into account \eqref{sf1}, we obtain
$$
\|u^{\delta,i}_{\eps,a}\|_{r, \partial_S \Om} \le C''
$$
for some $C''$ independent of $\delta$.
Since $g_\eps \in W^{1,1}([0,T];W^{1,p}(\Om;\R^2) \cap L^q(\Om;\R^2))$,
we obtain that for all $\tau \in [0,T]$
as $\delta \to 0$
$$
v_{\eps,a}^\delta(\tau) \to 0 \text{ strongly in }L^p(\Om;\msd),
$$
$$
w_{\eps,a}^\delta(\tau) \to 0 \text{ strongly in } L^{q}(\Om;\R^2),
$$
$$
z_{\eps,a}^\delta(\tau) \to 0 \text{ strongly in }
L^{r}(\partial_S \Om;\R^2).
$$
Moreover $\gamma^{\delta,\eps,a}_\fs \to 0$ and 
$\gamma^{\delta,\eps,a}_\gs \to 0$ as
$\delta \to 0$. Finally, by
\cite[Lemma 4.9]{DMFT}, as $\delta \to 0$
we have that for all $\tau \in [0,T]$
$$
|\langle \partial
\ws(\nabla u^\delta_{\eps,a}(\tau)+v_{\eps,a}^\delta(\tau)),
\nabla\dot{g}_\eps(\tau) \rangle
-\langle \partial
\ws(\nabla u^\delta_{\eps,a}(\tau)),
\nabla\dot{g}_\eps(\tau) \rangle| \to 0,
$$
$$
|\langle \partial
\fs(\tau)(u^\delta_{\eps,a}(\tau)+w_{\eps,a}^\delta(\tau)),
\dot{g}_\eps(\tau) \rangle
-\langle \partial \fs(\tau)(u^\delta_{\eps,a}(\tau)),
\dot{g}_\eps(\tau) \rangle| \to 0,
$$
$$
|\langle \partial
\gs(\tau)(u^\delta_{\eps,a}(\tau)+z_{\eps,a}^\delta(\tau)),
 \dot{g}_\eps(\tau) \rangle
-\langle \partial \gs(\tau)(u^\delta_{\eps,a}(\tau)),
\dot{g}_\eps(\tau) \rangle| \to 0,
$$
uniformly in $\eps,a$.
By the Dominated Convergence Theorem, we conclude that
$e^\delta_{\eps,a} \to 0$ as $\delta \to 0$ uniformly in $\eps$ and $a$,
and the proof is finished.
\end{proof}

\section{The approximation result}
\label{approx}
In this section we study the asymptotic behavior of the
discrete evolution obtained in Section
\ref{devol}.
Let us consider a given initial crack
$\Gamma^0 \in \gadm$ where $\gadm$ is defined as in
\eqref{gammaom}, and a boundary deformation
$g \in W^{1,1}\big([0,T], W^{1,p}(\Om;\R^2) 
\cap L^q(\Om;\R^2)\big)$.
Let $\Gamma^0_{\eps,a} \in \gadmea$ be an approximation
of $\Gamma^0$ in the sense of Proposition \ref{gammazero},
and let us consider
$$
g_\eps \in W^{1,1}\big([0,T], W^{1,p}(\Om;\R^2)
\cap L^q(\Om;\R^2)\big),
$$
such that
$$
g_\eps(t) \in \afeom \text{ for all }t \in [0,T],
$$
and such that for $\eps \to 0$
\begin{equation*}
g_\eps \to g
\quad \mbox{ strongly in }W^{1,1}([0,T],
W^{1,p}(\Om;\R^2) \cap L^q(\Om;\R^2)).
\end{equation*}
Let
\begin{equation*}
\{(u^{\delta,i}_{\eps,a},\Gamma^{\delta,i}_{\eps,a}),\,i=0, \ldots,N_\delta\}
\end{equation*}
be the discrete evolution relative to the initial crack
$\Gamma^0_{\eps,a}$ and
boundary data $g_\eps$ given by Proposition \ref{discrevol}.
We make the following piecewise constant interpolation in time:
\begin{equation}
\label{discrinterp}
u^\delta_{\eps,a}(t):=u^{\delta,i}_{\eps,a},
\quad\quad
\Gamma^\delta_{\eps,a}(t):=\Gamma^{\delta,i}_{\eps,a},
\quad\quad
g^\delta_\eps(t):=g_\eps(t^\delta_i)
\quad\quad \text{for } t^\delta_i \le t <t^\delta_{i+1},
\end{equation}
$i=0,\dots,N_\delta-1$,
and $u^\delta_{\eps,a}(T):=u^{\delta,N_\delta}_{\eps,a}$,
$\Gamma^\delta_{\eps,a}(T):=\Gamma^{\delta,N_\delta}_{\eps,a}$,
$g^\delta_\eps(T):=g_\eps(T)$.
\par
By Proposition \ref{dener}, for all $v \in \afeaomb$ we have
\begin{equation*}
\Eb(0)(u^{\delta}_{\eps,a}(0))
+\Es \left( \Sg{g^\delta_\eps(0)}{u^{\delta}_{\eps,a}(0)}
\setminus \Gamma^0_{\eps,a} \right)
\le \Eb(0)(v)+
\Es\left( \Sg{g^\delta_\eps(0)}{v}
\setminus \Gamma^0_{\eps,a}\right),
\end{equation*}
and for all $t \in [t^\delta_i, t^\delta_{i+1}[$ and for all
$v \in \afeaomb$
\begin{equation}
\label{piecemin2*}
\Eb(t^\delta_i)(u^{\delta}_{\eps,a}(t)) \le \Eb(t^\delta_i)(v) +
\Es \left( \Sg{g^\delta_\eps(t)}{v}
\setminus \Gamma^{\delta}_{\eps,a}(t) \right).
\end{equation}
Here $\Eb$ and $\Es$ are defined in
\eqref{elener} and \eqref{crackener}
respectively. Finally for all $0 \le s \le t \le 1$ we have
\begin{align}
\label{discrenergy2*}
\Eub(t^\delta_i)(u^{\delta}_{\eps,a}(t),\Gamma^\delta_{\eps,a}(t))
\le
&\Eub(s^\delta_j)(u^{\delta}_{\eps,a}(s),\Gamma^\delta_{\eps,a}(s))
+\int_{s^\delta_j}^{t^\delta_i}
\langle \partial \ws(\nabla u^\delta_{\eps,a}(\tau)),
\nabla \dot{g}_\eps(\tau) \rangle \,d\tau \\
\nonumber
&-\int_{s^\delta_j}^{t^\delta_i}
\dot{\fs}(\tau)(u^\delta_{\eps,a}(\tau))
-\int_{s^\delta_j}^{t^\delta_i}
\langle \partial \fs(\tau)(u^\delta_{\eps,a}(\tau)),
 \dot{g}_\eps(\tau) \rangle \,d\tau \\
\nonumber
&-\int_{s^\delta_j}^{t^\delta_i}
\dot{\gs}(\tau)(u^\delta_{\eps,a}(\tau))
-\int_{s^\delta_j}^{t^\delta_i}
\langle \partial \gs(\tau)(u^\delta_{\eps,a}(\tau)),
\dot{g}_\eps(\tau) \rangle \,d\tau
+e^\delta_{\eps,a},
\end{align}
where $s^\delta_j \le s <s^\delta_{j+1}$ and
$t^\delta_i \le t <t^\delta_{i+1}$,
$e^{\delta}_{\eps,a}$ is defined as in
\eqref{opiccolo}, and
$\Eub(t)(u,\Gamma)$ is as in \eqref{totalener}.
Recall that $e^{\delta}_{\eps,a} \to 0$ as $\delta \to 0$
uniformly in $\eps,a$.
\par
Comparing $u^{\delta}_{\eps,a}(t)$ with
$g^\delta_\eps(t)$ by \eqref{piecemin2*},
and in view of \eqref{wsb}, \eqref{fs2}, \eqref{fs3},
\eqref{gbou2}, \eqref{gbou3},
\eqref{MSdiscr} and \eqref{ik3},
by \eqref{discrenergy2*} with $s=0$
we deduce that there exists
$C' \in ]0,+\infty[$ such that for all
$t$, $\delta$, $\eps$ and $a$
\begin{equation}
\label{discrenergybis}
\|\nabla u^\delta_{\eps,a}(t)\|_p+\|u^\delta_{\eps,a}(t)\|_q+
\hs^1(\Gamma^\delta_{\eps,a}(t)) \le C'.
\end{equation}
By the time dependence of $\Eb(\cdot,\cdot)$,
in view of \eqref{discrenergybis}, by \eqref{piecemin2*} and 
\eqref{discrenergy2*} we have that there exists
$o^\delta_{\eps,a} \to 0$ as $\delta, \eps \to 0$
uniformly in $a$ such that for all $t \in [0,T]$ and
for all $v \in \afeaomb$
\begin{equation}
\label{piecemin2}
\Eb(t)(u^{\delta}_{\eps,a}(t)) \le \Eb(t)(v) +
\Es \left( \Sg{g^\delta_\eps(t)}{v}
\setminus \Gamma^{\delta}_{\eps,a}(t) \right)+o^\delta_{\eps,a},
\end{equation}
and for all $0 \le s \le t \le T$
\begin{align}
\label{discrenergy2}
\Eub(t)(u^\delta_{\eps,a}(t),\Gamma^\delta_{\eps,a}(t)) \le&
\Eub(s)(u^\delta_{\eps,a}(s),\Gamma^\delta_{\eps,a}(s))
+\int_{s}^{t}
\langle \partial \ws(\nabla u^\delta_{\eps,a}(\tau)),
\nabla \dot{g}_\eps(\tau) \rangle \,d\tau \\
\nonumber
&-\int_{s}^{t}
\dot{\fs}(\tau)(u^\delta_{\eps,a}(\tau))
-\int_{s}^{t}
\langle \partial \fs(\tau)(u^\delta_{\eps,a}(\tau)),
 \dot{g}_\eps(\tau) \rangle \,d\tau \\
\nonumber
&-\int_{s}^{t}
\dot{\gs}(\tau)(u^\delta_{\eps,a}(\tau))
-\int_{s}^{t}
\langle \partial \gs(\tau)(u^\delta_{\eps,a}(\tau)),
\dot{g}_\eps(\tau) \rangle \,d\tau
+o^\delta_{\eps,a}.
\end{align}
\par
Inequality \eqref{discrenergybis} gives a natural precompactness
of $(u^\delta_{\eps,a}(t))$ in
$GSBV^p_q(\Om;\R^2)$. The main result of the paper
is the following.

\begin{theorem}
\label{mainthm}
Let $\delta>0$, $\eps>0$, $a \in ]0,\frac{1}{2}[$, and
let $\left \{t \to \left( u_{\eps,a}^\delta(t),
\Gamma_{\eps,a}^\delta(t) \right)\,:\,t \in [0,T] \right\}$
be the discrete evolution
given by \eqref{discrinterp} relative to the initial crack
$\Gamma^0_{\eps,a}$ and the boundary data $g_\eps$.
Then there exist a quasistatic evolution
$\{t \to (u(t),\Gamma(t))\}$ in the sense of Theorem
\ref{qse} and sequences $\delta_n \to 0$, $\eps_n \to 0$,
$a_n \to 0$,
such that setting
$u_n(t):=u^{\delta_n}_{\eps_n,a_n}(t)$ and
$\Gamma_n(t):=\Gamma^{\delta_n}_{\eps_n,a_n}(t)$,
for all $t \in [0,T]$ the following facts hold.
\begin{itemize}
\item[]
\item[(a)]
For every $t \in [0,T]$, $(u_n(t))_{n \in \N}$
is weakly precompact in
$GSBV^p_q(\Om;\R^2)$, and every
accumulation point $\tilde{u}(t)$ is such that
$\Sg{g(t)}{\tilde{u}(t)} \tsub \Gamma(t)$,
\begin{equation}
\label{minimality}
\Eb(t)(\tilde{u}(t)) \le \Eb(t)(v)+
\Es \left(\Sg{g(t)}{v}\setminus \Gamma(t) \right)
\end{equation}
for all $v \in GSBV^p_q(\Om;\R^2)$ with
$S(v) \tsub \OmBb$, and
$$
\Eb(t)(u_n(t)) \to \Eb(t)(\tilde{u}(t)).
$$
Moreover there exists a subsequence of
$(\delta_n, \eps_n, a_n)_{n \in \N}$ (depending on $t$)
such that
$$
u_n(t) \weak u(t)
\quad
\text{ weakly in } GSBV^p_q(\Om; \R^2).
$$
\item[]
\item[(b)]
For every $t \in [0,T]$ we have
\begin{equation}
\label{convenergy}
\Eub(t)(u_n(t),\Gamma_n(t)) \to \Eub(t)(u(t),\Gamma(t));
\end{equation}
more precisely elastic and surface energies
converge separately, that is
\begin{equation}
\label{separateconv}
\Eb(t)(u_n(t)) \to \Eb(t)(u(t)),
\quad\quad
\Es(\Gamma_n(t)) \to \Es(\Gamma(t)).
\end{equation}
\end{itemize}
\end{theorem}

For the proof of Theorem \ref{mainthm} we need two
preliminary steps.
First of all, we fix $a$ and study the asymptotic
for $\delta,\eps \to 0$ (Lemma \ref{agamma}),
and then we let $a \to 0$ using a diagonal argument
(Lemma \ref{ato0}).

\begin{lemma}
\label{agamma}
Let $a$ be fixed, $t \in [0,T]$, and let $\delta_n \to 0$
and $\eps_n \to 0$.
There exists $\Gamma_a(t) \in \radm$ and a subsequence
of $(\delta_n,\eps_n)_{n \in \N}$ (which we denote with the same symbol),
such that the following facts hold:
\begin{itemize}
\item[(a)] if $w_n \in \afenab$ is such that
$\Sg{g^{\delta_n}_{\eps_n}(t)}{w_n}
\subseteq \Gamma^{\delta_n}_{\eps_n,a}(t)$
and
$$
w_n \weak w
\quad\quad
\text{weakly in }GSBV^p_q(\Om;\R^2),
$$
then we have
$$
\Sg{g(t)}{w} \tsub \Gamma_a(t);
$$
\item[]
\item[(b)]
there exists $\mu(a)$ with $\mu(a) \to 1$ as $a \to 0$ such that
for every accumulation point $u_a(t)$ of
$(u^{\delta_n}_{\eps_n,a}(t))_{n \in \N}$ for the weak convergence in
$GSBV^p_q(\Om;\R^2)$ and for all $v \in GSBV^p_q(\Om;\R^2)$
with $S(v) \tsub \OmBb$, we have
\begin{equation}
\label{aminimality}
\Eb(t)(u_a(t)) \le \Eb(t)(v)
+\mu(a) \Es \left( \Sg{g(t)}{v}
\setminus \Gamma_a(t) \right);
\end{equation}
moreover
\begin{equation}
\label{abulkconv}
\lim_n \Eb(t)(u^{\delta_n}_{\eps_n,a}(t))=
\Eb(t)(u_a(t));
\end{equation}
\item[]
\item[(c)] we have
$$
\Es(\Gamma_a(t)) \le \liminf_n \Es(\Gamma^{\delta_n}_{\eps_n,a}(t)).
$$
\end{itemize}
\end{lemma}

\begin{proof}
We now perform a variant of \cite[Theorem 4.7]{DMFT}.
Let $(\varphi_k)_{k \in \N}
\subseteq L^1(\Om;\R^2)$
be dense in $L^{1}(\Om;\R^2)$.
For every $\varphi_k$ and for every $m \in \N$, let
$v^{n,a}_{k,m}(t)$ be a minimum of
the problem
$$
\min\{
\|\nabla v\|_p+\|v\|_q+
m\|v-\varphi_k\|_1\,:\, v \in V^n_a\},
$$
where
$$
V^n_a:=\{v \in \afenab,\;
\Sg{g^{\delta_n}_{\eps_n}(t)}{v}
\subseteq \Gamma^{\delta_n}_{\eps_n,a}(t)\}.
$$
Since by \eqref{discrenergybis} we have
$\hs^1(\Gamma^{\delta_n}_{\eps_n,a}(t)) \le C'$,
by Theorem \ref{GSBVcompact}
there exists a subsequence of $(\delta_n,\eps_n)_{n \in \N}$
(which we denote with the same symbol)
such that $v^{n,a}_{k,m}(t)$ weakly converges to some
$v^a_{k,m}(t) \in GSBV^p_q(\Om;\R^2)$
as $n \to +\infty$ for all $k,m \in \N$. We set
\begin{equation}
\label{firstdefgammaat}
\Gamma_a(t):=
\bigcup_{k,m} \Sg{g(t)}{v^a_{k,m}(t)}.
\end{equation}
Let us see that $\Gamma_a(t)$ satisfies all the
properties of the lemma.
Clearly $\Gamma_a(t) \in \radm$ and point $(c)$
is a consequence of Theorem \ref{GSBVcompact}.
In particular by
\eqref{discrenergybis} we have that
\begin{equation}
\label{h1finitea}
\hs^1(\Gamma_a(t)) \le C'.
\end{equation}
\par
Let us come to point $(a)$. Let $w_n \in \afenab$ be such that
$\Sg{g^{\delta_n}_{\eps_n}(t)}{w_n}
\subseteq \Gamma^{\delta_n}_{\eps_n,a}(t)$
and
$$
w_n \weak w
\quad\quad
\text{weakly in }GSBV^p_q(\Om;\R^2).
$$
We claim that there exists $k_m \to +\infty$ such that
\begin{equation}
\label{convkm}
v^a_{k_m,m}(t) \weak w
\quad\quad
\text{weakly in }GSBV^p_q(\Om;\R^2).
\end{equation}
Then since
$\Sg{g(t)}{v^a_{k_m,m}(t)} \subseteq \Gamma_a(t)$ for all $m$
and in view of \eqref{h1finitea},
we deduce that
$\Sg{g(t)}{w} \tsub \Gamma_a(t)$. Let us prove \eqref{convkm}.
Fixed $m \in \N$, let us choose $k_m$ in such a way that
$$
m\|w-\varphi_{k_m}\|_1 \to 0.
$$
By minimality of $v^{n,a}_{k_m,m}(t)$ we have for all $n$
\begin{equation*}
\|\nabla v^{n,a}_{k_m,m}(t)\|_p
+\|v^{n,a}_{k_m,m}(t)\|_q
+m\|v^{n,a}_{k_m,m}(t)-\varphi_{k_m}\|_1
\le
\|\nabla w_n\|_p
+\|w_n\|_q
+m\|w_n-\varphi_{k_m}\|_1.
\end{equation*}
Passing to the limit in $n$, by lower semicontinuity
we get for some $C \ge 0$
\begin{equation*}
\|\nabla v^a_{k_m,m}(t)\|_p
+\|v^a_{k_m,m}(t)\|_q
+m\|v^a_{k_m,m}(t)-\varphi_{k_m}\|_{1}
\le C+m\|w-\varphi_{k_m}\|_1.
\end{equation*}
We deduce for $m \to +\infty$
$$
\|v^a_{k_m,m}(t)-\varphi_{k_m}\|_1 \to 0,
$$
which together with $\|\varphi_{k_m}-w\|_1 \to 0$
implies that
$$
v^a_{k_m,m}(t) \to w
\quad\quad
\text{strongly in }L^{1}(\Om;\R^2).
$$
Since
$$
\|\nabla v^a_{k_m,m}(t)\|_p
+\|v^a_{k_m,m}(t)\|_q \le
C+m\|w-\varphi_{k_m}\|_1 \le C+1
$$
for $m$ large, we have that
$v^a_{k_m,m}(t) \weak w$ weakly in $GSBV^p_q(\Om;\R^2)$,
and this proves \eqref{convkm}.
\par
Finally, let us come to point $(b)$.
Let $v \in GSBV^p_q(\Om;\R^2)$ with
$S(v) \subseteq \OmBb$, and let us
fix $k_1,\dots,k_s$ and $m_1,\dots, m_r$ in $\N$.
By Proposition \ref{piecetransf2}, there exists
$v_n \in \afenab$ such that
$$
\lim_n \Eb(t)(v_n)=\Eb(t)(v)
$$
and
\begin{multline*}
\limsup_n
\Es \left( \Sg{g^{\delta_n}_{\eps_n}(t)}{v_n} \setminus
\Gamma^{\delta_n}_{\eps_n,a}(t) \right) \le
\limsup_n \Es \big( \Sg{g^{\delta_n}_{\eps_n}(t)}{v_n}
\setminus
\bigcup_{i \le s,\,j \le r}
S(v^{n,a}_{k_i,m_j}) \big) \\
\le \mu(a) \Es \big( \Sg{g(t)}{v} \setminus
\bigcup_{i \le s,\,j \le r} S(v^a_{k_i,m_j}) \big),
\end{multline*}
where $\mu(a) \to 1$ as $a \to 0$.
Since the $k_i$'s and the $m_j$'s are arbitrary,
we obtain that
\begin{equation}
\label{limsupcrack}
\limsup_n \Es \left(
\Sg{g^{\delta_n}_{\eps_n}(t)}{v_n}
\setminus \Gamma^{\delta_n}_{\eps_n,a}(t)
\right) \le
\mu(a) \Es \left( \Sg{g(t)}{v}
\setminus \Gamma_a(t) \right).
\end{equation}
Let us suppose that
$u^{\delta_n}_{\eps_n,a}(t) \weak u_a(t)$
weakly in $GSBV^p_q(\Om;\R^2)$
along a suitable subsequence which we indicate by the same symbol.
By the minimality property \eqref{piecemin2},
comparing $u^{\delta_n}_{\eps_n,a}(t)$ with $v_n$ we get
\begin{equation}
\label{comparing}
\Eb(t)(u^{\delta_n}_{\eps_n,a}(t)) \le
\Eb(t)(v_n)+
\Es \left( \Sg{g^{\delta_n}_{\eps_n}(t)}{v_n}
\setminus \Gamma^{\delta_n}_{\eps_n,a}(t) \right)+o_n,
\end{equation}
with $o_n \to 0$ as $n \to +\infty$.
Then we have
\begin{multline*}
\Eb(t)(u_a(t)) \le \liminf_n
\Eb(t)(u^{\delta_n}_{\eps_n,a}(t)) \\
\le \limsup_n
\left( \Eb(t)(v_n)+
\Es \left( \Sg{g^{\delta_n}_{\eps_n}(t)}{v_n}
\setminus \Gamma^{\delta_n}_{\eps_n,a}(t) \right)
\right) \\
\le \Eb(t)(v)+ \limsup_n \Es \left(
\Sg{g^{\delta_n}_{\eps_n}(t)}{v_n} \setminus
\Gamma^{\delta_n}_{\eps_n,a}(t) \right) \\
\le \Eb(t)(v)+\mu(a) \Es \left( \Sg{g(t)}{v}
\setminus \Gamma_a(t) \right),
\end{multline*}
that is \eqref{aminimality} holds.
Choosing $v=u_a(t)$, passing to the limsup in
\eqref{comparing}, and taking into account
\eqref{limsupcrack} we obtain that
$$
\limsup_n \Eb(t)(u^{\delta_n}_{\eps_n,a}(t)) \le
\Eb(t)(u_a(t)).
$$
Since by \eqref{aminimality} $\Eb(t)(u_a(t))$ is independent
of the accumulation point $u_a(t)$, we conclude that
\eqref{abulkconv} holds.
\end{proof}

\begin{remark}
\label{evola}
{\rm
Using Lemma \ref{agamma}, it is possible to construct
an increasing family $\{t \to \Gamma_a(t)\,:\,t \in [0,T]\}$ and
a subsequence of $(\delta_n,\eps_n)_{n \in \N}$ such that
points $(a)$, $(b)$ and $(c)$ of Lemma \ref{agamma} hold
for every $t \in [0,T]$. This evolution
$\{t \to \Gamma_a(t)\,:\,t \in [0,T]\}$
can be considered as an approximate quasistatic evolution,
in the sense that it satisfies irreversibility, but satisfies
static equilibrium and nondissipativity up to a small error
due to the fact that $a$ is kept fixed.
The presence of $\mu(a)$ in the minimality property
\eqref{aminimality} takes into account the anisotropy in the
approximation of the surface energy:
in fact, since $a$ is kept fixed, the adaptive
edges of the triangulations $\treaom$ cannot recover all the possible
directions. The nondissipativity condition
up to a small error can be obtained
using the minimality property \eqref{aminimality} and
following \cite[Theorem 3.13]{DMFT} (estimate from below of
the total energy).
\par
The construction of $\{t \to \Gamma_a(t)\,:\,t \in [0,T]\}$
is the following.
If $D \subseteq [0,T]$ is countable and dense, by Lemma
\ref{agamma} and using a diagonalization argument,
we can find a subsequence of $(\delta_n,\eps_n)_{n \in \N}$ and an
increasing family $\Gamma_a(t) \in \radm$, $t \in D$, such
that points $(a)$, $(b)$ and $(c)$ hold for every $t \in D$.
Let us set for every $t \in [0,T]$
$$
\Gamma^+_a(t):= \bigcap_{s \ge t, s \in D} \Gamma_a(s).
$$
Clearly $\{t \to \Gamma^+_a(t)\,:\,t \in[0,T]\}$ is increasing, in
the sense that $\Gamma_a(s) \tsub \Gamma_a(t)$ for all
$s \le t$. As a consequence, the set
$J$ of discontinuity points of $\hs^1(\Gamma^+_a(t))$
is at most countable.
We can extract a further subsequence of $(\delta_n,\eps_n)_{n \in \N}$
such that $\Gamma_a(t)$ is determined
also for all $t \in J$ (notice that
$\Gamma_a(t) \tsub \Gamma^+_a(t)$). For all
$t \in [0,T] \setminus (D \cup J)$ we set
$\Gamma_a(t):=\Gamma^+_a(t)$.
We have that $\Gamma_a(t) \in \radm$ and
$\{t \to \Gamma_a(t)\,:\,t \in [0,T]\}$ is increasing.
\par
For $t \in D \cup J$, $\Gamma_a(t)$ satisfies by construction
points $(a)$, $(b)$ and $(c)$ of Lemma \ref{agamma}.
Let us consider the case
$t \in [0,T] \setminus (D \cup J)$.
\par
Concerning point $(a)$, we have that
$\Sg{g(t)}{w} \tsub \Gamma_a(s)$
for every $s \in D \cap [t,T]$, so that
passing to the intersection we get
$\Sg{g(t)}{u_a(t)} \tsub \Gamma_a(t)$.
\par
As for point $(b)$,
considering $s \in D \cap [0,t[$, for every $v \in
GSBV^p_q(\Om;\R^2)$ with $S(v) \tsub \OmBb$,
we have there exists $v_n \in \afenab$ such that
$$
\lim_n \Eb(t)(v_n)=\Eb(t)(v),
$$
and
\begin{multline*}
\limsup_n
\Es \left( \Sg{g^{\delta_n}_{\eps_n}(t)}{v_n}
\setminus \Gamma^{\delta_n}_{\eps_n,a}(t) \right) \le
\limsup_n
\Es \left( \Sg{g^{\delta_n}_{\eps_n}(t)}{v_n}
\setminus \Gamma^{\delta_n}_{\eps_n,a}(s) \right) \\
\le \mu(a) \Es \left( \Sg{g(t)}{v}
\setminus \Gamma_a(s) \right).
\end{multline*}
Then by minimality property \eqref{piecemin2} and passing
to the limit in $n$ we have
$$
\Eb(t)(u(t)) \le \Eb(t)(v)+\mu(a) \Es \left( \Sg{g(t)}{v}
\setminus \Gamma_a(s) \right).
$$
Letting $s \to t$ we get
that \eqref{aminimality} holds. Reasoning as in
Lemma \ref{agamma}, we get that also \eqref{abulkconv} holds.
\par
Finally, coming to point $(c)$,
we have that for all $s \in D \cap [0,t[$
$$
\liminf_n \Es(\Gamma^{\delta_n}_{\eps_n,a}(t)) \ge
\liminf_n \Es(\Gamma^{\delta_n}_{\eps_n,a}(s)) \ge
\Es(\Gamma_a(s)),
$$
so that letting $s \nearrow t$, and recalling that
$t$ is a continuity point for 
$\Es(\Gamma^{\delta_n}_{\eps_n,a}(\cdot))$,
we obtain that the lower semicontinuity holds.
}
\end{remark}

We can now let $a \to 0$.

\begin{lemma}
\label{ato0}
There exist a map
$\{t \to \Gamma(t) \in \radm, t \in [0,T]\}$
and sequences $\delta_n \to 0$,
$\eps_n \to 0$, $a_n \to 0$ such that
the following facts hold:
\begin{itemize}
\item[]
\item[(a)] $\Gamma^0 \tsub \Gamma(s) \tsub \Gamma(t)$ 
for all $0 \le s \le t \le T$;
\item[]
\item[(b)] for all $t \in [0,T]$, if $w_n \in \afenab$
with $\Sg{g^{\delta_n}_{\eps_n}(t)}{w_n} \subseteq
\Gamma^{\delta_n}_{\eps_n,a_n}(t)$ is such that
$$
w_n \weak w
\quad\quad
\text{weakly in }GSBV^p_q(\Om;\R^2),
$$
then we have
$$
\Sg{g(t)}{w} \tsub \Gamma(t);
$$
\item[]
\item[(c)] for all $t \in [0,T]$ and
for every accumulation point $u(t)$ of
$(u^{\delta_n}_{\eps_n,a_n}(t))_{n \in \N}$
for the weak convergence in $GSBV^p_q(\Om;\R^2)$
and for all $v \in GSBV^p_q(\Om;\R^2)$ with
$S(v) \tsub \OmBb$, we have
\begin{equation}
\label{minimality2}
\Eb(t)(u(t)) \le \Eb(t)(v)
+\Es \left( \Sg{g(t)}{v} \setminus \Gamma(t) \right),
\end{equation}
and
\begin{equation}
\label{bulkconv}
\Eb(t)(u(t))=
\lim_n
\Eb(t)(u^{\delta_n}_{\eps_n,a_n}(t));
\end{equation}
\item[]
\item[(d)] for all $t \in [0,T]$ we have
\begin{equation}
\label{liminfcrack}
\Es(\Gamma(t)) \le \liminf_n
\Es(\Gamma^{\delta_n}_{\eps_n,a_n}(t)).
\end{equation}
\end{itemize}
\end{lemma}

\begin{proof}
Let us consider $\delta_h \to 0$ and $\eps_h \to 0$.
Given $a \in ]0,\frac{1}{2}[$ and $t \in [0,T]$,
let $\Gamma_a(t)$ be the rectifiable set
given by Lemma \ref{agamma}.
Recall that by \eqref{firstdefgammaat}
we have
$$
\Gamma_a(t)= \bigcup_{k,m} \Sg{g(t)}{v^a_{k,m}(t)},
$$
where
$v^a_{k,m}(t)$ is the weak limit in $GSBV^p_q(\Om;\R^2)$
along a suitable subsequence depending on $a$
of a minimum $v^{h,a}_{k,m}(t)$ of
the problem
\begin{equation}
\label{defvhakm}
\min \{ \|\nabla v\|_p
+\|v\|_q
+m\|v-\varphi_k\|_1\,:\, v \in V^h_a(t)\},
\end{equation}
where $(\varphi_k)_{k \in \N} \subseteq L^1(\Om;\R^2)$
is dense in $L^1(\Om;\R^2)$ and
$$
V^h_a(t):=\{v \in \afenabp,
\Sg{g^{\delta_h}_{\eps_h}(t)}{v}
\subseteq \Gamma^{\delta_h}_{\eps_h,a}(t)\}.
$$
\par
Let $a_n \to 0$, and let $D:=\{t_j\,:\,j \in \N\}
\subseteq [0,T]$
be countable and dense with $0 \in D$.
Using a diagonal argument,
up to a subsequence of $(\delta_h,\eps_h)_{h \in \N}$,
we may suppose that for all $t \in D$ and for all $n$
$$
v^{h,a_n}_{k,m}(t) \weak v^{a_n}_{k,m}(t)
\quad\quad
\text{weakly in }GSBV^p_q(\Om;\R^2).
$$
Moreover, we may assume that for all $t \in D$ and for all $n$
$$
u^{\delta_h}_{\eps_h,a_n}(t) \weak u_{a_n}(t)
\quad\quad
\text{weakly in }GSBV^p_q(\Om;\R^2)
$$
with
$$
\Eb(t)(u^{\delta_h}_{\eps_h,a_n}(t))
\to \Eb(t)(u_{a_n}(t)).
$$
By Lemma \ref{agamma}, we have that $u_{a_n}(t)$
satisfies the minimality property
\eqref{aminimality}.
\par
Up to a subsequence of $(a_n)_{n \in \N}$, we may suppose that
for all $k,m$ and $t \in D$ we have
\begin{equation}
\label{convvkm}
v^{a_n}_{k,m}(t) \weak v_{k,m}(t)
\quad\quad
\text{weakly in }GSBV^p_q(\Om;\R^2),
\end{equation}
and
\begin{equation}
\label{convuan}
u_{a_n}(t) \weak u(t)
\quad\quad
\text{weakly in }GSBV^p_q(\Om;\R^2).
\end{equation}
For all $t \in D$, let us set
\begin{equation}
\label{defgammat1}
\Gamma(t):= \bigcup_{k,m} \Sg{g(t)}{v_{k,m}(t)}.
\end{equation}
By Proposition \ref{piecetransf2}, in view of the
minimality property \eqref{aminimality} and taking into account
that $\mu(a_n) \to 1$, we have that for all
$v \in GSBV^p_q(\Om;\R^2)$ with $S(v) \tsub
\OmBb$
\begin{equation}
\label{minaftera}
\Eb(t)(u(t)) \le \Eb(t)(v)+\Es( \Sg{g(t)}{v} \setminus \Gamma(t)),
\end{equation}
and as a consequence, we obtain
\begin{equation*}
\Eb(t)(u_{a_n}(t)) \to \Eb(t)(u(t)).
\end{equation*}
\par
We now perform the following diagonal argument. Choose
$\delta_{h_0},\eps_{h_0}$ in such a way that
\begin{multline*}
\|v^{h_0,a_0}_{0,0}(t_0)
-v^{a_0}_{0,0}(t_0)\|_1 +
\|u^{\delta_{h_0}}_{\eps_{h_0},a_0}(t_0)-u_{a_0}(t_0)
\|_1 \\
+|\Eb(t_0)(u^{\delta_{h_0}}_{\eps_{h_0},a_0}(t_0))-
\Eb(t_0)(u_{a_0}(t_0))|
\le 1.
\end{multline*}
Supposing to have constructed
$\delta_{h_n},\eps_{h_n}$, we choose
$\delta_{h_{n+1}},\eps_{h_{n+1}}$
in such a way that for all
$k \le n+1$, $m \le n+1$ and for all
$t_i$ with $1 \le i \le n+1$ we have
\begin{multline*}
\|v^{h_{n+1},a_{n+1}}_{k,m}(t_i)
-v^{a_{n+1}}_{k,m}(t_i)\|_1 +
\|u^{\delta_{h_{n+1}}}_{\eps_{h_{n+1}},a_{n+1}}(t_i)
-u_{a_{n+1}}(t_i)\|_1 \\
+|\Eb(t_i)(u^{\delta_{h_{n+1}}}_{\eps_{h_{n+1}},a_{n+1}}(t_i))-
\Eb(t_i)(u_{a_{n+1}}(t_i))|
\le \frac{1}{n+1}.
\end{multline*}
Let us set $\delta_n:=\delta_{h_n}$ and $\eps_n:=\eps_{h_n}$, and
let us prove that $\Gamma(t)$ defined in \eqref{defgammat1}
satisfies the properties of the Lemma.
We have immediately that $\Gamma(t) \in \radm$.
\par
Concerning point $(d)$, notice that
$$
\Gamma^{\delta_n}_{\eps_n,a_n}(t)=
\bigcup_{m,k} \Sg{g^{\delta_n}_{\eps_n}(t)}{v^{h_n,a_n}_{k,m}(t)},
\quad\quad
\Gamma(t)=\bigcup_{m,k} \Sg{g(t)}{v_{k,m}(t)},
$$
and that for all $k,m$
$$
v^{h_n,a_n}_{k,m}(t) \weak v_{k,m}(t)
\quad\quad
\text{weakly in }GSBV^p_q(\Om;\R^2);
$$
then \eqref{liminfcrack} is a consequence
of Theorem \ref{GSBVcompact}. In particular,
by \eqref{discrenergybis}, we get that
\begin{equation}
\label{h1finite}
\hs^1(\Gamma(t)) \le C'.
\end{equation}
\par
Let us come to point $(b)$. Let $w_n \in \afenab$
with $\Sg{g^{\delta_n}_{\eps_n}(t)}{w_n}
\subseteq \Gamma^{\delta_n}_{\eps_n,a_n}(t)$
be such that
$$
w_n \weak w
\quad\quad
\text{weakly in }GSBV^p_q(\Om;\R^2).
$$
For every $m \in \N$, let us choose $k_m$ in such a way that
$$
m\|w-\varphi_{k_m}\|_1 \to 0.
$$
By minimality of $v^{h_n,a_n}_{k_m,m}(t)$ we have for all $n$
\begin{equation*}
\|\nabla v^{h_n,a_n}_{k_m,m}(t)\|_p
+\|v^{h_n,a_n}_{k_m,m}(t)\|_q
+m\|v^{h_n,a_n}_{k_m,m}(t)-\varphi_{k_m}\|_1
\le
\|\nabla w_n\|_p
+\|w_n\|_q
+m\|w_n-\varphi_{k_m}\|_1.
\end{equation*}
By construction of $h_n$, and in view of \eqref{convvkm}, we have
$$
v^{h_n,a_n}_{k_m,m}(t) \weak v_{k_m,m}(t)
\quad\quad
\text{weakly in }GSBV^p_q(\Om;\R^2).
$$
Then passing to the limit in $n$, by lower semicontinuity
we get for some $C \ge 0$
\begin{equation*}
\|\nabla v_{k_m,m}(t)\|_p
+\|v_{k_m,m}(t)\|_q
+m\|v_{k_m,m}(t)-\varphi_{k_m}\|_1 \le
C+m\|w-\varphi_{k_m}\|_1.
\end{equation*}
We deduce for $m \to +\infty$
$$
\|v_{k_m,m}(t)-\varphi_{k_m}\|_1 \to 0,
$$
which together with
$\|\varphi_{k_m}-w\|_1 \to 0$
implies that
$$
v_{k_m,m}(t) \to w \quad
\text{ strongly in }L^{1}(\Om;\R^2).
$$
Since
$$
\|\nabla v_{k_m,m}(t)\|_p
+\|v_{k_m,m}(t)\|_q
\le C+m\|w-\varphi_{k_m}\|_1 \le C+1
$$
for $m$ large, we have that
$$
v_{k_m,m}(t) \weak w
\quad\quad
\text{weakly in }GSBV^p_q(\Om;\R^2).
$$
Since $\Sg{g(t)}{v_{k_m,m}(t)} \subseteq \Gamma(t)$
for all $m$, and
since $\hs^1(\Gamma(t))<C'$, we deduce that
$\Sg{g(t)}{w} \tsub \Gamma(t)$.
\par
Coming to point $(c)$, we have that \eqref{bulkconv}
holds by construction.
Moreover \eqref{minimality2} holds in view 
of \eqref{minaftera} and by the fact that
$u^{\delta_n}_{\eps_n,a_n}(t)$ weakly
converges in $GSBV^p_q(\Om;\R^2)$ to $u(t)$
defined in \eqref{convuan}.
\par
In order to prove point $(a)$, notice that if $s \le t$
with $s,t \in D$, we have for all $k,m \in \N$
that
$$
\Sg{g^{\delta_n}_{\eps_n}(t)}{v^{h_n,a_n}_{k,m}(s)+
g^{\delta_n}_{\eps_n}(t)-g^{\delta_n}_{\eps_n}(s)}
\subseteq \Gamma^{\delta_n}_{\eps_n,a_n}(s)
\subseteq \Gamma^{\delta_n}_{\eps_n,a_n}(t),
$$
and
$$
v^{h_n,a_n}_{k,m}(s)+
g^{\delta_n}_{\eps_n}(t)-g^{\delta_n}_{\eps_n}(s)
\weak
v_{k,m}(s)+g(t)-g(s)
\quad\quad
\text{weakly in }GSBV^p_q(\Om;\R^2),
$$
where $v^{h,a}_{k,m}(s)$ and $v_{k,m}(s)$
are defined in \eqref{defvhakm} and
\eqref{convvkm}.
By point $(b)$ we deduce that
$$
\Sg{g(t)}{v_{k,m}(s)+g(t)-g(s)}=\Sg{g(s)}{v_{k,m}(s)}
\tsub \Gamma(t).
$$ 
Then by the definition of $\Gamma(s)$
we get $\Gamma(s) \tsub \Gamma(t)$. Finally,
by the same argument, we deduce $\Gamma^0 \tsub \Gamma(s)$.
\par
In order to deal with all $t \in [0,T]$,
we proceed as in Remark \ref{evola}.
For all $t \in [0,T] \setminus D$ let us set
$$
\Gamma^+(t):= \bigcap_{s \ge t, s \in D} \Gamma(s).
$$
Clearly $\Gamma^+(t) \in \radm$ and satisfies point $(a)$,
so that the set $J$ of discontinuity points
of $\hs^1(\Gamma^+(\cdot))$ is at most countable.
We can then extract a further subsequence
of $(\delta_n,\eps_n, a_n)_{n \in \N}$
such that $\Gamma(t)$ is determined
also for all $t \in J \setminus D$
(notice that $\Gamma(t) \tsub \Gamma^+(t)$).
For all $t \in [0,T] \setminus (D \cup J)$ we set
$\Gamma(t):=\Gamma^+(t)$.
We have that $\Gamma(t) \in \radm$ and
that $\Gamma(t)$ satisfies point $(a)$.
Let us see that $\Gamma(t)$ satisfies also
points $(b)$, $(c)$ and $(d)$ also for $t \in
[0,T] \setminus (D \cup J)$.
\par
Concerning point $(b)$, for every
accumulation point  $u(t)$
of $(u^{\delta_n}_{\eps_n,a_n}(t))_{n \in \N}$
for the weak convergence in
$GSBV^p_q(\Om;\R^2)$, by the first part of the proof,
we have that $\Sg{g(t)}{u(t)} \tsub \Gamma(s)$
for all $s \in D$ with
$s \ge t$, so that passing to the intersection,
we get that $\Sg{g(t)}{u(t)} \tsub \Gamma(t)$.
\par
Let us come to point $(c)$.
Let
$$
u_j(t):=u^{\delta_{n_j}}_{\eps_{n_j},a_{n_j}}(t) \weak u(t)
\quad\quad
\text{weakly in }GSBV^p_q(\Om;\R^2)
$$
along a subsequence $n_j \nearrow +\infty$. Let us set
$\Gamma_j:=\Gamma^{\delta_{n_j}}_{\eps_{n_j},a_{n_j}}$ and
$g_j:=g^{\delta_{n_j}}_{\eps_{n_j}}$.
Up to a further subsequence there exists
$s_j \in D$ with $s_j \nearrow t$, and such that setting
$u_j(s_j):=u^{\delta_{n_j}}_{\eps_{n_j},a_{n_j}}(s_j)$, we have
\begin{equation}
\label{controlconv}
\|u_j(s_j)-u(s_j)\|_1
+
|\Eb(s_j)(u_j(s_j))-\Eb(s_j)(u(s_j))| \to 0.
\end{equation}
We have that there exists
$u^*(t) \in GSBV^p_q(\Om;\R^2)$ such that up
to a subsequence
$$
u(s_j) \weak u^*(t)
\quad\quad
\text{weakly in }GSBV^p_q(\Om;\R^2).
$$
By the minimality property \eqref{minimality2} of $u(s_j)$,
for all $v \in GSBV^p_q(\Om;\R^2)$
with $S(v) \tsub \Omb_B$, we have that
$$
\Eb(s_j)(u(s_j)) \le \Eb(s_j)(v-g(t)+g(s_j))
+\Es(\Sg{g(t)}{v} \setminus \Gamma(s_j)).
$$
Passing to the limit in $j$ we have that
for all $v \in GSBV^p_q(\Om;\R^2)$
with $S(v) \tsub \Omb_B$
\begin{equation}
\label{min*}
\Eb(t)(u^*(t)) \le \Eb(t)(v)+\Es(\Sg{g(t)}{v}
\setminus \Gamma(t)).
\end{equation}
As a consequence of the stability of this unilateral
minimality property, it follows that
\begin{equation*}
\Eb(s_j)(u(s_j)) \to \Eb(t)(u^*(t)).
\end{equation*}
By \eqref{controlconv} we get
$$
u_j(s_j) \weak u^*(t)
\quad\quad
\text{weakly in }GSBV^p_q(\Om;\R^2),
$$
and
\begin{equation}
\label{conv*}
\Eb(s_j)(u_j(s_j)) \to \Eb(t)(u^*(t)).
\end{equation}
By \eqref{piecemin2}, comparing $u_j(t)$
with $u_j(s_j)-g_j(s_j)+g_j(t)$,
taking into account that
$$
\Sg{g_j(s_j)}{u_j(s_j)}
\subseteq
\Gamma_j(s_j) \subseteq \Gamma_j(t),
$$
we obtain
\begin{equation*}
\Eb(t)(u_j(t)) \le
\Eb(s_j)(u_j(s_j))+o_j
\end{equation*}
where $o_j \to 0$ as $j \to +\infty$.
Passing to the limit in $j$ we have by \eqref{conv*}
$$
\Eb(t)(u(t)) \le \liminf_j \Eb(t)(u_j(t)) \le
\limsup_j \Eb(t)(u_j(t)) \le \Eb(t)(u^*(t)).
$$
By \eqref{min*} we deduce that
\eqref{minimality2} holds. Moreover
we have that $\Eb(t)(u(t))=\Eb(t)(u^*(t))$
and that $\Eb(t)(u(t))$ is independent of the
accumulation point $u(t)$. Then we deduce that
\eqref{bulkconv} holds.
\par
Finally, concerning point $(d)$, we have that for all
$s \in D \cap [0,t[$
$$
\liminf_n \Es(\Gamma^{\delta_n}_{\eps_n,a_n}(t)) \ge
\liminf_n \Es(\Gamma^{\delta_n}_{\eps_n,a_n}(s)) \ge
\Es(\Gamma(s)),
$$
so that letting $s \nearrow t$ we obtain \eqref{liminfcrack}.
The proof is now complete.
\end{proof}

We can now prove Theorem \ref{mainthm}.

\begin{proof}[\underline{\rm PROOF OF THEOREM \ref{mainthm}}]
Let $(\delta_n,\eps_n,a_n)_{n \in \N}$ and
$\{t \to \Gamma(t) \in \radm, t \in [0,T]\}$
be given by Lemma \ref{ato0}. For all $t \in [0,T]$, let us set
$$
u_n(t):=u^{\delta_n}_{\eps_n,a_n}(t), \quad\quad
\Gamma_n(t):= \Gamma^{\delta_n}_{\eps_n,a_n}(t).
$$
Let us see that it is possible to choose an accumulation point
$u(t) \in GSBV^p_q(\Om;\R^2)$ of $(u_n(t))_{n \in \N}$
such that $\{t \to (u(t),\Gamma(t))\,:\,
t \in [0,T]\}$ is a quasistatic growth of brittle fractures
in the sense of Dal Maso-Francfort-Toader.
Let us set
\begin{multline*}
\vartheta_n(s):=
\langle \partial \ws(\nabla u_n(s)),
\nabla \dot{g}_{\eps_n}(s) \rangle \\
-\dot{\fs}(s)(u_n(s))-
\langle \partial \fs(s)(u_n(s)),
\dot{g}_{\eps_n}(s) \rangle \\
-\dot{\gs}(s)(u_n(s))-
\langle \partial \gs(s)(u_n(s)),
\dot{g}_{\eps_n}(s) \rangle.
\end{multline*}
By growth conditions of $\ws, \fs, \gs$ and
by \eqref{discrenergybis} we have that
there exists $\psi \in L^1(0,T)$ such that
$\vartheta_n(s) \le \psi(s)$ for all $n$.
Let us consider
$$
\vartheta(s):= \limsup_n \vartheta_n(s).
$$
By \cite[Theorem 5.5 and Lemma 4.11]{DMFT},
for every $s \in [0,T]$
there exists $u(s)$ accumulation point of
$(u_n(s))_{n \in \N}$ for the weak convergence
in $GSBV^p_q(\Om;\R^2)$
such that
\begin{multline*}
\vartheta(s):=
\langle \partial \ws(\nabla u(s)),
\nabla \dot{g}(s) \rangle \\
-\dot{\fs}(s)(u(s))-
\langle \partial \fs(s)(u(s)),
\dot{g}(s) \rangle \\
-\dot{\gs}(s)(u(s))-
\langle \partial \gs(s)(u(s)),
\dot{g}(s) \rangle.
\end{multline*}
Applying  Fatou's Lemma (in the limsup version) to
\eqref{discrenergy2} with $s=0$, we have that
$$
\Eub(t)(u(t),\Gamma(t))
\le \limsup_n \Eub(0)(u_n(0),\Gamma_n(0))+
\int_0^t \vartheta(s) \,ds.
$$
By Proposition \ref{gammazero}, we have that
$\limsup_n \Eub(0)(u_n(0),\Gamma_n(0))=\Eub(0)(u(0),\Gamma(0))$, so
that we get
$$
\Eub(t)(u(t),\Gamma(t)) \le
\Eub(0)(u(0),\Gamma(0))+ \int_0^t \vartheta(s) \,ds.
$$
Moreover, again by \cite[Theorem 3.13]{DMFT},
$$
\Eub(t)(u(t),\Gamma(t))
\ge \Eub(0)(u(0),\Gamma(0))
+\int_0^t \vartheta(s) \,ds,
$$
so that
\begin{equation}
\label{nondiss}
\Eub(t)(u(t),\Gamma(t))=
\Eub(0)(u(0),\Gamma(0))+
\int_0^t \vartheta(s) \,ds.
\end{equation}
We deduce that $\{t \to (u(t),\Gamma(t))\,:\,t \in [0,T]\}$
is a quasistatic growth of brittle fractures:
in fact by Lemma \ref{ato0} we get that
$\Gamma(\cdot)$ is increasing, and for
$t \in [0,T]$ $(u(t),\Gamma(t)) \in AD(g(t))$
and the static equilibrium holds; moreover
the nondissipativity condition is given by \eqref{nondiss}.
\par
Let us see that points $(a)$ and $(b)$ of Theorem \ref{mainthm}
holds. By \eqref{discrenergybis}, $(u_n(t))_{n \in \N}$
is weakly precompact in $GSBV^p_q(\Om;\R^2)$
for all $t \in [0,T]$.
Moreover by Lemma \ref{ato0} every accumulation point
$\tilde{u}(t)$ of $(u_n(t))_{n \in \N}$ for
the weak convergence in
$GSBV^p_q(\Om;\R^2)$ is such that
$\Sg{g(t)}{\tilde{u}(t)} \subseteq \Gamma(t)$
and the minimality property \eqref{minimality} holds.
Moreover we have
$$
\Eb(t)(\tilde{u}(t))=
\lim_n \Eb(t)(u_n(t)).
$$
Since $\Eb(t)(\tilde{u}(t))$
is independent of the particular accumulation point
$\tilde{u}(t)$, we have that point $(a)$ is proved.
\par
Let us come to point $(b)$.
Taking into account \eqref{bulkconv} and
\eqref{liminfcrack}, for all $t \in [0,T]$ we have
$$
E(t) \le \liminf_n E_n(t) \le
\limsup_n E_n(t) \le E(0)+ \int_0^t \vartheta(s) \,ds=E(t),
$$
so that \eqref{convenergy} holds.
Moreover we deduce that separate convergence of elastic and
surface energies holds at any time, so that
\eqref{separateconv} is proved.
The proof is now concluded.
\end{proof}


\section{The strictly convex case}
\label{convexcase}
In this section we assume that the function $W(x,\xi)$
is strictly convex in $\xi$ for a.e. $x \in \Om$
and that the function $F(t,x,z)$ is strictly
convex in $z$ for all $t \in [0,T]$ and for a.e.
$x \in \Om$: as a consequence, the elastic energy
$\Eb(t,v)$ is strictly
convex in $v$ for all $t \in [0,T]$, and a stronger
approximation result is available.

\begin{theorem}
\label{main2}
Let $g \in W^{1,1}\big([0,T], W^{1,p}(\Om;\R^2)
\cap L^q(\Om;\R^2)\big)$ and let
$$
g_\eps \in W^{1,1}\big([0,T], W^{1,p}(\Om;\R^2)
\cap L^q(\Om;\R^2)\big),
\quad\quad
g_\eps(t) \in \afeom
\quad
\text{for all }t \in [0,T]
$$
be such that for $\eps \to 0$
$$
g_\eps \to g
\quad\quad
\text{strongly in }W^{1,1}([0,T], W^{1,p}(\Om;\R^2)
\cap L^q(\Om;\R^2)).
$$
Let $\Gamma^0 \in {\mathbf \Gamma}(\Om)$ be an initial crack
and let
$\Gamma^0_{\eps,a}$ be its approximation in the sense of
Proposition \ref{gammazero}.
Let us suppose that
\begin{align}
\label{assumptions}
W(x,\cdot) &\mbox{ is strictly convex for a.e. }x \in \Om, \\
\nonumber
F(t,x,\cdot) &\mbox{ is strictly convex for a.e. }
(t,x) \in [0,T] \times \Om.
\end{align}
Given  $\delta>0$, $\eps>0$, $a \in ]0,\frac{1}{2}[$,
let $\left \{t \to \left( u_{\eps,a}^\delta(t),
\Gamma_{\eps,a}^\delta(t) \right)\,:\,t \in [0,T] \right\}$
be the piecewise constant interpolation of the discrete
evolution given by Proposition \ref{discrevol}
relative to the initial crack $\Gamma^0_{\eps,a}$
and the boundary data $g_{\eps}$.
Then there exists a quasistatic evolution
$\{t \to (u(t),\Gamma(t))\,:\,t \in [0,T]\}$ relative
to the initial
crack $\Gamma^0$ and the boundary data $g$ in the
sense of Theorem \ref{qse}, and
sequences $\delta_n \to 0$, $\eps_n \to 0$, $a_n \to 0$,
such that setting
$$
u_n(t):=u^{\delta_n}_{\eps_n,a_n}(t),
\quad\quad
\Gamma_n(t):=\Gamma^{\delta_n}_{\eps_n,a_n}(t),
$$
for all $t \in [0,T]$ the following facts hold:
\begin{itemize}
\item[]
\item[(a)] $\nabla u_n(t) \to \nabla u(t)$ strongly
in $L^p(\Om;\msd)$ and
$u_n(t) \to u(t)$ strongly in
$L^q(\Om;\R^2)$;
\item[]
\item[(b)]
$\Eub(t)(u_n(t),\Gamma_n(t)) \to \Eub(t)(u(t),\Gamma(t))$,
and in particular elastic and
surface energies converge separately, that is
\begin{equation*}
\Eb(t)(u_n(t)) \to \Eb(t)(u(t)),
\quad \quad \Es(\Gamma_n(t)) \to \Es(\Gamma(t)).
\end{equation*}
\end{itemize}
\end{theorem}

\begin{proof}
Let us consider the sequence
$(\delta_n, \eps_n, a_n)_{n \in \N}$
and the quasistatic growth of brittle fractures
$\{t \to (u(t),\Gamma(t))\,:\,t \in [0,T]\}$
given in Theorem \ref{mainthm}. Under assumptions
\eqref{assumptions},
we have that $u(t)$ is uniquely determined, because
by \eqref{minimality} $u(t)$ minimizes
$$
\min \{ \Eb(t)(v) \,:\, v \in GSBV^p_q(\Om;\R^2),\;
\Sg{g(t)}{v} \subseteq \Gamma(t) \},
$$
and $\Eb(t)(\cdot)$ is strictly convex.
We conclude by point $(a)$ of
Theorem \ref{mainthm} that $u_n(t) \weak u(t)$ weakly in
$GSBV^p_q(\Om;\R^2)$.
Point $(b)$ is a direct consequence of
Theorem \ref{mainthm}. By the convergence of the
elastic energy, we deduce that
\begin{align*}
\lim_n \int_{\Om} W(x, \nabla u_n(t)) \,dx&=
\int_{\Om} W(x, \nabla u(t)) \,dx, \\
\lim_n \int_\Om F(t,x,u_n(t))\,dx&=
\int_\Om F(t,x,u(t))\,dx.
\end{align*}
By the assumption on the strict convexity
of $W$ and $F$ we deduce by \cite{Bre2}
$$
\nabla u_n(t) \to \nabla u(t)
\quad\quad
\text{strongly in }L^p(\Om;\msd),
$$
and
$$
u_n(t) \to u(t)
\quad\quad
\text{strongly in }L^q(\Om;\R^2).
$$
Point $(a)$ is now proved, and the proof is concluded.
\end{proof}

\bigskip
\bigskip
\centerline{ACKNOWLEDGMENTS}
\bigskip\noindent
The authors wish to thank Gianni Dal Maso
for many helpful and interesting discussions.


\end{document}